\documentclass[a4paper,12pt]{article}
\usepackage[T2A]{fontenc}
\usepackage[cp1251]{inputenc}
\begin{document}

\title{Homotopism of homological complexes over nonassociative
algebras with metagroup relations.}

\author{S.V. Ludkowski.}

\date{12 December 2020}
\maketitle
\begin{abstract}
The article is devoted to homological complexes. Smashly graded
modules and complexes are studied over nonassociative algebras with
metagroup relations. Smashed tensor products of homological
complexes are investigated. Their homotopisms and homologisms are
scrutinized.
\par \textbf{Keywords:} nonassociative algebra;
smashly graded; homology; complex; homotopism; homologism;
cohomology; metagroup
\end{abstract}

\section{Introduction.}
\par Nonassociative algebras and algebras with group relations are very important in different branches of
mathematics and its applications (see, for example,
\cite{cartaneilenbb56,gotogrossb,shangy2012,shangy2015,shangy2017}
and references therein). For example, octonions and generalized
Cayley-Dickson algebras played very important role in mathematics
and quantum field theory
\cite{allcja98,baez,dickson,kansol,schaeferb}. They are used not
only in algebra and noncommutative geometry, but also in
noncommutative analysis and PDEs, particle physics, mathematical
physics, mathematical analysis, operator theory and their
applications in natural sciences including physics and quantum field
theory (see
\cite{baez,dickson,frenludkfejms18}-\cite{guertzeb,kansol}-\cite{ludkcvee13}
and references therein). \par A multiplicative law of their
canonical generators is nonassociative and leads to a more general
notion of a metagroup instead of a group \cite{ludlmla18} (see also
Appendix). They were used in \cite{ludlmla18} for investigations of
automorphisms and derivations of nonassociative algebras. A
remarkable  fact was proved in the 20-th century that the nontrivial
geometry exists if and only if there exists the corresponding unital
quasigroup (see \cite{bruckb,kiechlb,pickert,smithb} and references
therein).
\par An extensive area of investigations of
PDEs intersects with cohomologies and deformed cohomologies
\cite{pommb}. Therefore, it is important to develop this area over
octonions, Cayley-Dickson algebras and more general metagroup
algebras. This means that metagroup algebras include as the
particular case the Cayley-Dickson algebras. This article is devoted
to nonassociative algebras with metagroup relations. It is worth to
notice that a class of metagroups differs substantially from a class
of groups. Indeed, a metagroup may be nonassociative or power
non-associative or nonalternative. Moreover, left or right inverse
elements in the metagroup may not exist or it may contain elements
for which left and right inverse elements do not coincide.
\par On the other hand, algebras are frequently studied using
cohomology theory. But the already developed cohomology theory
operates with associative algebras. It was investigated by
Hochschild and other authors
\cite{bourbalgbhomol,bredonb67,cartaneilenbb56,hochschild46}, but it
is not applicable to nonassociative algebras. In some particular
cases of nonassociative algebras such as Lie algebras, pre-Lie
algebras, flexible algebras, alternative algebras homology theory
was developed for needs of studies their structure (see, for
example, \cite{chapliver2001,dzhzus2011,gotogrossb,remmgoze} and
references therein). It is necessary to note that classes of these
algebras are quite different from classes of generalized
Cayley-Dickson algebras and nonassociative algebras with metagroup
relations.
\par Previously cohomologies of loop spaces on quaternion and octonion manifolds
were studied in \cite{ludwrgrijmgta}. They have specific features in
comparison with the case of complex manifolds. This is especially
caused by the noncommutativity of the quaternion skew field and the
nonassociativity of the octonion algebra. Principles of cohomology
theory for nonassociative algebras with metagroup relations were
described in \cite{ludkctnaax19}.
\par This article is devoted to homological complexes. Smashly graded
modules and complexes are studied over nonassociative algebras with
metagroup relations. Smashed tensor products of homological
complexes are investigated. Their homotopisms and homologisms are
scrutinized.
\par All main results of this paper are obtained for the first time.

\section{Tensor products of complexes for nonassociative
algebras with metagroup relations.}
\par {\bf  Definition 2.1.} Let $\mathcal T$ be an associative unital
ring, $G$ be a metagroup and $A={\mathcal T}[G]$ be a metagroup
algebra. Let also $B$ be a unital right $A$-module (see Definitions
A1 and A2 in Appendix) such that $B$ has also a structure of a
two-sided ${\mathcal T}$-module. Suppose that $B$ has a
decomposition $B=\sum_{g\in G}B_g$ as a two-sided ${\mathcal
T}$-module, where $B_g$ is a two-sided ${\mathcal T}$-module, and
satisfies the following conditions:
\par $(2.1.1)$ $B_gh=B_{gh}$;
\par $(2.1.2)$ $(bh)x_g=b(hx_g)$ and $x_g(bh)=(x_gh)b$ and $bx_g=x_gb$;
\par $(2.1.3)$ $(x_gh)s={\sf t}_3(g,h,s) x_g(hs)$;
\par $(2.1.4)$ $(bc)x=b(cx)$, $(bx)c=b(xc)$, $(xb)c=x(bc)$
\\ for every $h, g, s$
in $G$, $b$ and $c$ in $\mathcal T$, $x_g\in B_g$ and $x\in B$. \par
Then a right $A$-module $B$ will be called $G$-graded. Similarly are
considered left and two-sided modules, $(A,\mbox{}^1A)$-bimodules.
For $(A,A)$-bimodules it can be also shortly written $A$-bimodule or
two-sided $A$-module. Suppose that $B$ is an $A$-bimodule and
\par $(2.1.5)$ there exists a ${\mathcal T}$-bilinear mapping $B\times B\ni (x,y)\mapsto
xy\in B$ such that
\par $x(y+z)=xy+xz$ and $(y+z)x=yx+zx$ and $(bx)y=b(xy)$ and
$(xb)y=x(by)$ and $(xy)b=x(yb)$\\ for all $x$, $y$, $z$ in $B$,
$b\in \mathcal T$;
\par $(2.1.6)$ $(x_gy_h)z_s={\sf t}_3(g,h,s) x_g(y_hz_s)$ and
$x_gy_h\in B_{gh}$
\\ for every $g$, $h$, $s$
in $G$, $x_g\in B_g$, $y_h\in B_h$, $z_s\in B_s$. \par  Then we call
$B$ a $G$-graded algebra over $A$ (or a $G$-graded $A$-algebra). The
algebra $B$ is called unital if and only if
\par $(2.1.7)$ $B$ has a unit element $1=1_B$ such that $1_Bx=x$ and
$x1_B=x$ for each $x\in B$.

\par Assume that $B$ is the $A$-algebra and $P\subseteq B$, where
$A={\mathcal T}[G]$ is the metagroup algebra. We put
\par $(2.1.8)$ $Com _B(P) := \{ x\in B: ~ \forall b\in P, ~ xb=bx \} $;
\par $(2.1.9)$ $N_{B,l}(P) := \{ x\in B: ~ \forall b\in P, ~ \forall c\in P, ~ (xb)c=x(bc) \}
$;
\par $(2.1.10)$ $N_{B,m}(P) := \{ x\in B: ~ \forall b\in P, ~ \forall c\in P, ~ (bx)c=b(xc)
\} $;
\par $(2.1.11)$ $N_{B,r}(P) := \{ x\in B: ~ \forall b\in P, ~ \forall c\in P, ~ (bc)x=b(cx)
\} $;
\par $(2.1.12)$ $N_B(P) := N_{B,l}(P)\cap N_{B,m}(P)\cap N_{B,r}(P)$;
\par $(2.1.13)$ ${\sf C}_B(P) := Com _B(P)\cap N_B(P)$.
\par Then $Com _B(P)$, $~N_B(P)$ and ${\sf C}_B(P)$ are called a commutant,
a nucleus and a centralizer correspondingly of the algebra $B$
relative to a subset $P$ in $B$. Instead of $Com_B(B)$, $N_B(B)$ or
${\sf C}_B(B)$ it will be also written shortly  $Com(B)$, $N(B)$ or
${\sf C}(B)$ correspondingly. We put $P_e=P\cap B_e$ and $P_{\sf
C}:=P\cap B_{\sf C}$ for the $G$-graded $A$-algebra $B$, where
$B_{\sf C}=\sum_{g\in {\sf C}(G)} B_g$.

\par {\bf Lemma 2.2.} {\it Let $B$ be a $G$-graded $A$-algebra (see
Definition 2.1). Then $N(B)=\sum_{g\in N(G)} B_g$ and $N(B)$ is an
associative ${\mathcal T}$-algebra.}
\par {\bf Proof.} We put $Y=\sum_{g\in N(G)} B_g$. From Definition
2.1 it follows that $B$ has also a structure of a ${\mathcal
T}$-algebra and $Y$ is a ${\mathcal T}$-subalgebra in $B$. From
Conditions $(2.1.6)$ and $(2.1.12)$ we infer that $Y$ is the
associative ${\mathcal T}$-algebra and $Y\subset N(B)$.
\par On the other, hand Conditions $(2.1.9)$-$(2.1.12)$ imply tat
$N(B)$ is the associative ${\mathcal T}$-subalgebra in $B$. From
Definition 2.1 we deduce that \par $N(B)=\sum_{g\in G}[B_g\cap
N(B)]$.
\par Assume that $x\in N(B)$, $x\ne 0$. Then $x=\sum_{g\in
P(x)}x_g$,  where $P(x)\subset G$, $x_g\in B_g\cap N(B)$, $x_g\ne 0$
for each $g\in P(X)$. If $h\in (G\setminus N(G))\cap P(x)$, $y_h\in
B_h$, $y_h\ne 0$, then Conditions $(2.1.1)$ and $(2.1.6)$ imply that
$y_h\notin N(B)$. Therefore, $P(x)\subset N(G)$, consequently,
$N(B)\subset Y$. Thus $Y=N(B)$.

\par {\bf Definition 2.3.} Suppose that $\mbox{}^jB$ is a unital
$\mbox{}^jG$-graded $\mbox{}^jA$-algebra, where
$\mbox{}^jA={\mathcal T}[\mbox{}^jG]$ is a metagroup algebra for
each $j\in \{ 1, 2,...  \} $, ${\mathcal T}$ is an associative
unital ring. Suppose that $X$ is a $\mbox{}^1G$-graded left
$\mbox{}^1B$-module and $Y$ is a $\mbox{}^2G$-graded left
$\mbox{}^2B$-module (see also Definition A3).
\par  Suppose also that $f: X\to Y$ is a map such that $f$ is
\par $(2.3.1)$ a left ${\mathcal T}$-homomorphism and
$f(X_g)\subseteq Y_{f_{\bf \iota }(g)}$ for each $g\in \mbox{}^1G$,
where $f_{\bf \iota }: \mbox{}^1G\to \mbox{}^2G$ is a homomorphism
of metagroups:
\par $(2.3.2)$ $f_{\bf \iota }(gh)=f_{\bf \iota }(g)f_{\bf \iota }(h)$ and $f_{\bf \iota }(g\setminus
h)=f_{\bf \iota }(g)\setminus f_{\bf \iota }(h)$ and \par $f_{\bf
\iota }(g/h)=f_{\bf \iota }(g)/f_{\bf \iota }(h)$ for every $g$ and
$h$ in $G$.
\par The map $f: X\to Y$ satisfying conditions $(2.3.1)$ and
$(2.3.2)$ will be called a $(\mbox{}^1G,\mbox{}^2G)$-graded left
${\mathcal T}$-homomorphism of the left modules $X$ and $Y$. If the
ring is specified, it may be shortly written homomorphism instead of
${\mathcal T}$-homomorphism. Symmetrically is defined a
$(\mbox{}^3G,\mbox{}^4G)$-graded right homomorphism of a right
$\mbox{}^3B$-module $X$ and a right $\mbox{}^4B$-module $Y$. For a
$(\mbox{}^1B,\mbox{}^3B)$-bimodule $X$ and a
$(\mbox{}^2B,\mbox{}^4B)$-bimodule $Y$ if a map $f: X\to Y$ is
$(\mbox{}^1G,\mbox{}^2G)$-graded left and
$(\mbox{}^3G,\mbox{}^4G)$-graded right homomorphism, then $f$ will
be called a
$((\mbox{}^1G,\mbox{}^2G),(\mbox{}^3G,\mbox{}^4G))$-graded
${\mathcal T}$-homomorphism of bimodules $X$ and $Y$.

\par Assume that $X$ is a left $\mbox{}^1B$-module and $Y$ is a left $\mbox{}^2B$-module
and $f: X\to Y$ is a map such that
\par $(2.3.3)$ $f: X\to Y$ is a left ${\mathcal T}$-homomorphism and
$f(ax+by)=f_{\bf \iota }(a)f(x)+f_{\bf \iota }(b)f(y)$, where
\par $(2.3.4)$ $f_{\bf \iota }$ is a ${\mathcal T}$-homomorphism from
$\mbox{}^1B$ into $\mbox{}^2B$ such that \par $f_{\bf \iota }:
\mbox{}^1G\to \mbox{}^2G$, $f_{\bf \iota }: N(\mbox{}^1B)\to
N(\mbox{}^2B)$, $f_{\bf \iota }: \mbox{}^1B\to \mbox{}^2B$ is
injective, $f_{\bf \iota }(1_{\mbox{}^1B})=1_{\mbox{}^2B}$ and
$f_{\bf \iota }: {\mathcal T}[N(\mbox{}^1G)]\to {\mathcal
T}[N(\mbox{}^2G)]$ and
\par $(2.3.5)$ $f_{\bf \iota }(ab)=f_{\bf \iota }(a)f_{\bf \iota }(b)$ and
$f_{\bf \iota }(g\setminus h)=f_{\bf \iota }(g)\setminus
f_{\bf \iota }(h)$ and \par $f_{\bf \iota }(g/h)=f_{\bf \iota
}(g)/f_{\bf \iota }(h)$ and $f_{\bf \iota }(pa+bs)=pf_{\bf \iota
}(a)+f_{\bf \iota }(b)s$
\par for every $g$ and $h$ in $\mbox{}^1G$, $ ~ a$ and $b$ in $\mbox{}^1B$, $~ p$
and $s$ in ${\mathcal T}$, where $\mbox{}^jA$ is embedded into
$\mbox{}^jB$ as $\mbox{}^jA1_{\mbox{}^jB}$, $1_{\mbox{}^jB}$ is a
unit element in the unital $\mbox{}^jA$-algebra $\mbox{}^jB$,
$\mbox{}^jG$ is naturally emebedded into $\mbox{}^jA$ as
$\mbox{}^jG1_{\mathcal T}$, $1_{\mathcal T}$ is the unit in
${\mathcal T}$.
\par If $f$ satisfies conditions $(2.3.3)$-$(2.3.5)$, then $f$
will be called a $(\mbox{}^1B,\mbox{}^2B)$-generic left homomorphism
of left modules $X$ and $Y$. For right modules a
$(\mbox{}^1B,\mbox{}^2B)$-generic right homomorphism is defined
analogously. If $X$ is a $(\mbox{}^1B,\mbox{}^3B)$-bimodule and $Y$
is a $(\mbox{}^2B,\mbox{}^4B)$-bimodule and $f$ is a
$(\mbox{}^1B,\mbox{}^2B)$-generic left and
$(\mbox{}^3B,\mbox{}^4B)$-generic right homomorphism, then $f$ will
be called a
$((\mbox{}^1B,\mbox{}^2B),(\mbox{}^3B,\mbox{}^4B))$-generic
homomorphism of bimodules $X$ and $Y$. \par An
$(\mbox{}^1B,\mbox{}^2B)$-generic left homomorphism $f$ of left
modules $X$ and $Y$ such that $f_{\bf \iota }$ is an epimorphism of
$\mbox{}^1B$ onto $\mbox{}^2B$  we will call
$(\mbox{}^1B,\mbox{}^2B)$-epigeneric.
\par If additionally the homomorphism $f_{\bf \iota }$ is surjective
in $(2.3.4)$ and $f_{\bf \iota }^{-1}: \mbox{}^2B\to \mbox{}^1B$ is
the ${\cal T}$-homomorphism, then $f_{\bf \iota }$ is called an
isomorphism of $\mbox{}^1B$ with $\mbox{}^2B$ (or automorphism if
$\mbox{}^1B=\mbox{}^2B$).
\par A $(\mbox{}^1B,\mbox{}^2B)$-generic left homomorphism $f$ of left modules $X$
and $Y$ such that $f_{\bf \iota }$ is the isomorphism of
$\mbox{}^1B$ with $\mbox{}^2B$ will be called
$(\mbox{}^1B,\mbox{}^2B)$-exact.
\par In particular, if $\mbox{}^1G=\mbox{}^2G$, then
"$(\mbox{}^1G,\mbox{}^1G)$-graded" or
"$(\mbox{}^1B,\mbox{}^1B)$-generic" will shortened to
"$\mbox{}^1G$-graded" or "$\mbox{}^1B$-generic" correspondingly,
etc. If $f_{\bf \iota }$ is an automorphism of $\mbox{}^1G$ (or of
$\mbox{}^1B$ correspondingly), then the $\mbox{}^1G$-graded (or the
$\mbox{}^1B$-generic) left ${\mathcal T}$-homomorphism from $X$ into
$Y$ will be called $\mbox{}^1G$-exact (or $\mbox{}^1B$-exact
correspondingly). Similarly $\mbox{}^3G$-exact or $\mbox{}^3B$-exact
right homomorphisms of right modules and
$(\mbox{}^1G,\mbox{}^3G)$-exact or $(\mbox{}^1B,\mbox{}^3B)$-exact
homomorphisms of bimodules are defined (as shortening of
$((\mbox{}^1G,\mbox{}^1G),(\mbox{}^3G,\mbox{}^3G))$-exact or of
$((\mbox{}^1B,\mbox{}^1B),(\mbox{}^3B,\mbox{}^3B))$-exact
homomorphisms of bimodules).

\par If $X$ and $Y$ are $G$-graded $B$-algebras and
$f$ is a $G$-graded (or $G$-exact or $B$-generic or $B$-exact)
homomorphism from $X$ into $Y$ considered as $B$-bimodules  and in
addition the following condition is satisfied
\par $(2.3.6)$ $f(vx)=f(v)f(x)$ for each $x$ and $v$ in $X$, \\ then
$f$ will be called a $G$-graded (or $G$-exact or $B$-generic or
$B$-exact correspondingly) homomorphism of the $B$-algebras, where
$B=\mbox{}^jB$ for some $j\in \{ 1, 2,... \} $.

\par {\bf Definition 2.4.} Assume that $G$ is a metagroup, $A={\mathcal T}[G]$ is
a metagroup algebra $A={\mathcal T}[G]$ and $X$ is a two-sided
$B$-module, where $\mathcal T$ is an associative unital ring, $B$ is
a unital $G$-graded $A$-algebra. We denote by $G^n$ the $n$-fold
direct product of $G$ with itself such that $G^n$ is a metagroup,
where $n\ge 1$ is a natural number. We consider a two-sided
$N(B)$-module $X_{ \{ g_1,...,g_n \} _{q(n)}}$ for each
$g_1$,...,$g_n$ in $G$ and a vector $q(n)$ indicating an order of
pairwise multiplications in the braces $ \{ g_1,...,g_n \} $ (see
Definition 1 in \cite{ludkctnaax19}). Suppose that $X$ has the
following decomposition \par $(2.4.1)$ $X=\sum_{g_1\in G,..., g_n\in
G} X_{ \{ g_1,...,g_n \} _{l(n)}}$
\\ as the two-sided $N(B)$-module, where $ \{ g_1 \} _{l(1)}=g_1$,
$ ~ \{ g_1, g_2 \} _{l(2)}=g_1g_2$ and by induction $ \{
g_1,...,g_n, g_{n+1} \} _{l(n+1)}= \{ g_1,...,g_n \} _{l(n)}g_{n+1}$
for each $n\ge 1$. Assume also that $X$ satisfies the following
conditions:
\par $(2.4.2)$  there exists a left $N(B)$-linear and a right $N(B)$-linear isomorphism \par
$\theta (g_1,...,g_n;q(n),v(n)) : X_{ \{ g_1,...,g_n \} _{v(n)}}\to
X_{ \{ g_1,...,g_n \} _{q(n)}}$ such that
\par $\theta (g_1,...,g_n;q(n),v(n))(x_{ \{ g_1,...,g_n \}
_{v(n)}})= t_n(g_1,...,g_n;q(n),v(n)) x_{ \{ g_1,...,g_n \}
_{v(n)}}$ \\ for each $x_{ \{ g_1,...,g_n \} _{v(n)}}\in X_{ \{
g_1,...,g_n \} _{v(n)}}$, where $t_n(g_1,...,g_n;q(n),v(n))\in {\bf
\Psi }$ is such that \par $\{ g_1,...,g_n \} _{q(n)} =
t_n(g_1,...,g_n;q(n),v(n)) \{ g_1,...,g_n \} _{v(n)}$,
\par $t_n(g_1,...,g_n;q(n),v(n))=t_n(g_1,...,g_n;q(n),v(n)| id)$ \\ (see
also Lemma 1 and Example 2 in \cite{ludkctnaax19});
\par $(2.4.3)$ there exist left $N(B)$-linear and right $N(B)$-linear
isomorphisms \par
$\theta _l(g_0,g_1,...,g_n;l(n),l(n)) : g_0X_{ \{ g_1,...,g_n \}
_{l(n)}}\to X_{ \{ (g_0g_1),...,g_n \} _{l(n)}}$ and \par $\theta
_r(g_1,...,g_n,g_{n+1};l(n),l(n)) : X_{ \{ g_1,...,g_n \}
_{l(n)}}g_{n+1}\to X_{ \{ g_1,...,(g_ng_{n+1}) \} _{l(n)}}$
\par $(2.4.4)$
$(bg_0)x_{\{ g_1,...,g_n \} _{l(n)}}=b(g_0x_{\{ g_1,...,g_n \}
_{l(n)}})$ and \par $x_{\{ g_1,...,g_n \} _{l(n)}}(bg_{n+1})=(x_{\{
g_1,...,g_n \} _{l(n)}}g_{n+1})b$ and \par $bx_{\{ g_1,...,g_n \}
_{l(n)}}=x_{ \{ g_1,...,g_n \} _{l(n)}}b$,
\par $(2.4.5)$ $(g_0g_{n+1})x_{ \{ g_1,...,g_n \}
_{l(n)}}={\sf t}_3(g_0,g_{n+1},g) g_0(g_{n+1}x_{ \{ g_1,...,g_n \}
_{l(n)}})$ and \par $(g_0x_{ \{ g_1,...,g_n \} _{l(n)}})g_{n+1}=
{\sf t}_3(g_0,g,g_{n+1}) g_0(x_{ \{ g_1,...,g_n \} _{l(n)}}g_{n+1})$
and
\par $(x_{ \{ g_1,...,g_n \}
_{l(n)}}g_0)g_{n+1} = {\sf t}_3(g,g_0,g_{n+1}) x_{ \{ g_1,...,g_n \}
_{l(n)}}(g_0g_{n+1})$
\\ for every $b\in {\sf C}(B)$, $ ~ x_{ \{ g_1,...,g_n \} _{l(n)}}\in X_{ \{ g_1,...,g_n \} _{l(n)}}$,
elements $g_0$, $g_1$,...,$g_n$, $g_{n+1}$ in the metagroup $G$,
vectors $q(n)$ and $v(n)$ indicating orders of pairwise
multiplications, where $g= \{ g_1,...,g_n \} _{l(n)}$. Then a
two-sided $B$-module $X$ satisfying conditions $(2.4.1)$-$(2.4.5)$
will be called smashly $G^n$-graded. For short it also will be said
"$G^n$-graded" instead of "smashly $G^n$-graded". In particular, if
the module $X$ is $G^n$-graded and splits into a direct sum
\par $(2.4.6)$ $\quad X=\bigoplus_{(g_1,...,g_n)\in V}X_{ \{ g_1,...,g_n \}
_{l(n)}}$ \\ of two-sided $N(B)$-submodules $X_{ \{ g_1,...,g_n \}
_{l(n)}}$, where $V\subset G^n$, then we will say that that $X$ is
directly $G^n$-graded.
\par Similarly are defined $G^n$-graded left and right $B$-modules.
\par The $G^n$-graded left (or right or two-sided) $B$-module $X$ is called essentially
$G^n$-graded with $n\ge 1$, if for each $g_1$,...,$g_n$ in $G$ there
exists $h_n$ in $G$ such that $X_{ \{ g_{\sigma (1)},...,g_{\sigma
(n)} \} } \ne X_{ \{ h_{\sigma (1)},...,h_{\sigma (n)} \} }$ for
each $\sigma \in S_n$, where $h_j=g_j$ for each $j\le n-1$ if $n>1$;
$~ S_n$ is the symmetric group (i.e. of all permutations $\sigma $
of the finite set $ \{ 1,...,n \} $).

\par {\bf Lemma 2.5.} {\it Let $X$ be a smashly $G^n$-graded two-sided
(or left, or right) $B$-module, $1\le m \le n$. Then $X$ can be
supplied with a smashly $G^m$-graded two-sided (or left, or right
respectively) $B$-module structure.}
\par {\bf Proof.} The smashly $G^n$-graded two-sided $B$-module $X$
has the decomposition $(2.4.1)$ as the two-sided $N(B)$-module. By
virtue of Lemma 2.2 $N(B)$ is the associative ${\mathcal T}$-algebra
such that $N(B)=\sum_{g\in N(G)}B_g$ and $N(B)$ is contained in $B$.
The case $m=n$ is trivial. It remains the case $1\le m<n$. We put
\par $(2.5.1)$ $Y_{ \{ s, g_{n-m+2},...,g_n \} _{l(m)} } =\sum_{s= \{
g_1,...,g_{n-m+1} \} _{l(n-m+1)}; ~ g_1\in G,..., g_{n-m+1}\in G}
X_{ \{ g_1,...,g_n \} _{l(n)}}$
\\ for $2<m<n$
\par $(2.5.2)$ $Y_{ \{ s, g_n \} _{l(2)}} =\sum_{s= \{
g_1,...,g_{n-1} \} _{l(n-1)}; ~ g_1\in G,..., g_{n-1} \in G} X_{ \{
g_1,...,g_n \} _{l(n)}}$
\\ for $2=m<n$,
\par $(2.5.3)$ $Y_{ \{ s \} _{l(1)}} =\sum_{s= \{
g_1,...,g_{n} \} _{l(n)}; ~ g_1\in G,..., g_n\in G} X_{ \{
g_1,...,g_n \} _{l(n)}}$
\\ for $1=m<n$. Then we define
\par $(2.5.4)$ $Y=\sum_{g_1\in G,..., g_m\in G} Y_{ \{ g_1,...,g_m \}
_{l(m)}}$. \par Notice that if both vectors $q(n)$ and $s(n)$ in
$g_1,...,g_{n-m+1}$ terms correspond to $q(n-m+1)$ order of
multiplication, that is \par $(2.5.5)$ $\{ g_1,...,g_{n} \} _{q(n)}=
\{ \{ g_1,...,g_{n-m+1} \} _{q(n-m+1)}, g_{n-m+2},..., g_n \}
_{u(m)}$
\par $\{ g_1,...,g_{n} \} _{v(n)}= \{ \{
g_1,...,g_{n-m+1} \} _{q(n-m+1)}, g_{n-m+2},..., g_n \} _{w(m)}$
\\ with the corresponding vectors $u(m)$ and $w(m)$,
then \par $t_n(g_1,...,g_n;q(n),v(n))= t_m(p, g_{n-m+2},..., g_n
;u(m),w(m))$,\\ where $p = \{ g_1,...,g_{n-m+1} \} _{q(n-m+1)}$.
Therefore from $(2.4.1)$-$(2.4.5)$ and $(2.5.1)$-$(2.5.5)$ it
follows that $Y$ is a smashly $G^m$-graded two-sided $B$-module.
This construction supplies $X$ with the smashly $G^m$-graded
two-sided $B$-module structure.

\par {\bf Definition 2.6.} If $B$ is a $G^n$-graded two-sided
$A$-module, $1<n\in {\bf N}$, if also $B$ supplied with the
$G$-graded two-sided $A$-module structure by Lemma 2.5 is a
$G$-graded $A$-algebra (see Definition 2.1), then $B$ will be called
a $G^n$-graded $A$-algebra.

\par {\bf Lemma 2.7.} {\it Let $B$ be a smashly $G^k$-graded $A$-algebra,
let also $X$ be a smashly $G^m$-graded two-sided (or left, or right)
$B$-module, $1\le k\le m \le n\in {\bf N}$ (see Definitions 2.4 and
2.6). Then there exists a smashly $G^n$-graded two-sided (or left,
or right respectively) $B$-module $Y$ such that $Y$ can be supplied
with a smashly $G^m$-graded two-sided (or left, or right
respectively) $B$-module structure relative to which it is
isomorphic with $X$. Moreover, $X$ and $Y$ are isomorphic if
considered as two-sided (or left, or right respectively)
$N(B)$-modules.}
\par {\bf Proof.} The case of $m=n$ is trivial taking $X=Y$.
Let now $1\le m<n$. We choose $Y$ such that
\par $(2.7.1)$ $Y =\sum_{g_1\in G,..., g_n\in G} Y_{ \{
g_1,...,g_n \} _{l(n)}}$ with
\par $(2.7.2)$ $Y_{ \{
g_1,...,g_n \} _{l(n)}}= X_{ \{ h, g_{n-m+2},...,g_n \} _{l(m)} } $,
\\ for $2<m$, where $h=\{ g_1,...,g_{n-m+1} \} _{l(n-m+1)}$,
\par $(2.7.3)$ $Y_{ \{
g_1,...,g_n \} _{l(n)}}= X_{ \{ h, g_n \} _{l(2)} } $ \\
for $m=2$, where $h=\{ g_1,...,g_{n-1} \} _{l(n-1)}$,
\par $(2.7.4)$ $Y_{ \{
g_1,...,g_n \} _{l(n)}}=X_{ \{ h \} _{l(1)} }$\\ for $m=1$, where
$h=\{ g_1,...,g_{n} \} _{l(n)}$.
\par Notice that the sum in $(2.7.1)$ may generally may be not
direct (see also Definition 2.4). Formulas $(2.7.1)$-$(2.7.4)$ imply
that $X$ and $Y$ are isomorphic as two-sided (or left, or right
respectively) $N(B)$-modules. In view of Lemma 2.5 this $Y$ can be
supplied with a smashly $G^m$-graded two-sided (or left, or right
respectively) $B$-module structure. Relative to the latter structure
$Y$ is isomorphic with $X$ by $(2.7.2)$-$(2.7.4)$.
\par By virtue of Lemmas 2.2 and 2.5
$N(B)=\sum_{g\in N(G)} B_g$ and $N(B)$ is an associative ${\mathcal
T}$-algebra. As the ${\mathcal T}$-algebra $N(B)$ is the subalgebra
in $B$ considered as the ${\mathcal T}$-algebra. Then $X$ and $Y$
are the two-sided (or left, or right respectively) $N(B)$-modules,
since $X$ and $Y$ are the two-sided (or left, or right respectively)
$B$-modules. From Definitions 2.1 and 2.4 it follows that
$a(bx)=(ab)x$ for each $a$ and $b$ in $N(B)$ and $x\in X$, if $X$ is
the left $B$-module, $(xb)a=x(ba)$ if $X$ is the right $B$-module.
Therefore, $(2.7.1)$-$(2.7.4)$ imply that $X$ and $Y$ are isomorphic
if considered as the two-sided (or left, or right respectively)
$N(B)$-modules.

\par {\bf  Definition 2.8.} Assume that ${\mathcal C}$ is a $G$-graded
$B$-bimodule (see Definitions 2.1, 2.3), where $A={\mathcal T}[G]$
is the metagroup algebra, $B$ is a unital $G$-graded $A$-algebra,
${\mathcal T}$ is a commutative associative unital ring. Assume also
that ${\mathcal C}$ is ${\bf Z}$-graded with a gradation $({\mathcal
C}_n: n\in {\bf Z})$, where ${\bf Z}$ is the ring of all integers,
${\mathcal C}^n={\mathcal C}_{-n}$. Suppose also that ${\mathcal
C}^n$ is $G^{|n|+2}$-graded $B$-bimodule for each integer $n$ (see
Definition 2.4). Suppose also that $d: {\mathcal C}\to {\mathcal C}$
is a ${\bf Z}$-graded $B$-generic homomorphism of degree $-1$ such
that $dd=0$. Then $d$ is called a differential of a $G$-graded
$B$-complex $({\mathcal C},d)$ or ${\mathcal C}$. Shortly it may be
written a differential complex $({\mathcal C},d)$ or ${\mathcal C}$,
if $G$ and $B$ are specified.
\par Let $(\mbox{}^1{\mathcal C},\mbox{}^1d)$ be a $\mbox{}^1G$-graded $\mbox{}^1B$-complex.
If $\psi : ({\mathcal C},d)\to (\mbox{}^1{\mathcal C},\mbox{}^1d)$
is a homomorphism of degree $0$ such that it is
$((B,\mbox{}^1B),(B,\mbox{}^1B))$-generic (or
$((B,\mbox{}^1B),(B,\mbox{}^1B))$-exact) with $\mbox{}^1d\circ \psi
= \psi \circ d$, then $\psi $ is called a homomorphism (or a
$((B,\mbox{}^1B),(B,\mbox{}^1B))$-exact homomorphism
correspondingly) of complexes.
\par Similarly are defined complexes and their homomorphisms
in other cases: if ${\mathcal C}$ is a $G$-graded left $B$-module or
right $B$-module.
\par {\bf Remark 2.9.} In this work mainly essentially $G^n$-graded
$B$-modules are considered (see also Proposition 2.26 below). Lemmas
2.5 and 2.7 serve in order to encompass other cases for their
simultaneous treatment.
\par Henceforth differential $G$-graded $B$-complexes and unital
$G$-graded $A$-algebras $B$ are considered for metagroups $G$ and
metagroup algebras $A={\mathcal T}[G]$, if something other will not
be outlined. \par  Definition 2.8 means that $d_n: {\mathcal C}_n\to
{\mathcal C}_{n-1}$ and $d_nd_{n+1}=0$, also $d^n: {\mathcal C}^n\to
{\mathcal C}^{n+1}$ and $d^nd^{n-1}=0$ for each $n\in {\bf Z}$.
Examples of $G$-graded $A$-complexes are provided by
\cite{ludkctnaax19} (see also Proposition 1 and Theorem 1 there),
where $d_n$ and $d^n$ are left and right $A$-homomorphisms,
consequently, $A$-exact, hence $A$-generic.
\par Then ${\mathcal Z}({\mathcal C},d)=ker ~(d) $ is a module of cycles, ${\mathcal B}({\mathcal C},d)= Im ~(d)$
is a module of boundaries. They are $B$-bimodules $G$-graded and
${\bf Z}$-graded such that ${\mathcal Z}_n({\mathcal C},d)=ker ~
(d_n)$, $~ {\mathcal B}_n({\mathcal C},d)= Im ~(d_{n+1})$, $ ~
{\mathcal Z}_n({\mathcal C},d)={\mathcal Z}^{-n}({\mathcal C},d)$, $
~ {\mathcal B}_n({\mathcal C},d)={\mathcal B}^{-n}({\mathcal C},d)$
for each $n\in {\bf Z}$.
\par The homomorphism $\psi : ({\mathcal C},d)\to (\mbox{}^1{\mathcal C},\mbox{}^1d)$
of $G$-graded $B$- and $\mbox{}^1B$-complexes correspondingly means
that ${\mbox{}^1d}_n\circ \psi _n= \psi _{n-1}\circ d_n$ and
${\mbox{}^1d}^n\circ \psi ^n = \psi ^{n+1}\circ d_n$ for each $n\in
{\bf Z}$.
\par Assume that $({\mathcal C},d)$ is a differential $G$-graded $B$-complex,
$(\mbox{}^1{\mathcal C},\mbox{}^1d)$ is a differential
$\mbox{}^1G$-graded $\mbox{}^1B$-complex, $(\mbox{}^2{\mathcal
C},\mbox{}^2d)$ is a differential $\mbox{}^2G$-graded
$\mbox{}^2B$-complex.
\par Then $\psi ({\mathcal Z}({\mathcal C}))\subset {\mathcal Z}(\mbox{}^1{\mathcal C})$ and
$\psi ({\mathcal B}({\mathcal C}))\subset {\mathcal
B}(\mbox{}^1{\mathcal C})$. This induces homomorphisms ${\mathcal
Z}(\psi ): {\mathcal Z}({\mathcal C})\to {\mathcal
Z}(\mbox{}^1{\mathcal C})$ and ${\mathcal B}(\psi ): {\mathcal
B}({\mathcal C})\to {\mathcal B}(\mbox{}^1{\mathcal C})$ and $H(\psi
): H({\mathcal C})\to H(\mbox{}^1{\mathcal C})$. Their ${\bf
Z}$-homogeneous components are ${\mathcal Z}_n(\psi )$, ${\mathcal
Z}^n(\psi )$,...,$H^n(\psi )$. If $\psi $ and $\phi $ are
homomorphisms of $({\mathcal C},d)$ into $(\mbox{}^1{\mathcal
C},\mbox{}^1d)$, then $\psi +\phi $ is a homomorphism from
$({\mathcal C},d)$ into $(\mbox{}^1{\mathcal C},\mbox{}^1d)$ such
that ${\mathcal Z}(\psi +\phi )={\mathcal Z}(\psi ) + {\mathcal
Z}(\phi )$ and ${\mathcal B}(\psi +\phi )={\mathcal B}(\psi ) +
{\mathcal B}(\phi )$ and $H(\psi +\phi )=H(\psi ) + H(\phi )$.
Moreover, $b\psi $ is a homomorphism from $({\mathcal C},d)$ into
$(\mbox{}^1{\mathcal C},\mbox{}^1d)$ for each $b\in {\mathcal T}$
with ${\mathcal Z}(b\psi )=b{\mathcal Z}(\psi )$, ${\mathcal
B}(b\psi )=b{\mathcal Z}(\psi )$ and $H(b\psi )=bH(\psi )$.
\par If $H(\psi )$ is bijective, then it is called a
$((B,\mbox{}^1B),(B,\mbox{}^1B))$-generic (or
$((B,\mbox{}^1B),(B,\mbox{}^1B))$-exact) homologism respectively or
shortly homologism. If $H(\psi )=0$, then the differential
$G$-graded $B$-complex $({\mathcal C},d)$ is called null
homological. If $H_n({\mathcal C})=0$ ($H^n({\mathcal C})=0$), then
$({\mathcal C},d)$ is called acyclic descending degree $n$
(ascending degree $n$ respectively).
\par If $\chi : (\mbox{}^1{\mathcal C}, \mbox{}^1d)\to (\mbox{}^2{\mathcal C},\mbox{}^2d)$ is a homomorphisms, then $\chi
\circ \psi $ is a homomorphism from $({\mathcal C},d)$ into
$(\mbox{}^2{\mathcal C},\mbox{}^2d)$ with ${\mathcal Z}(\chi \circ
\psi )={\mathcal Z}(\chi )\circ {\mathcal Z}(\psi )$ and ${\mathcal
B}(\chi \circ \psi )={\mathcal B}(\chi )\circ {\mathcal B}(\psi )$
and $H(\chi \circ \psi )=H(\chi )\circ H(\psi )$.
\par Analogously are considered complexes, their homomorphisms and
homologisms in other cases: if ${\mathcal C}$ is a $G$-graded left
$B$-module or right $B$-module. Then Lemma 2.10, Definitions 2.11,
2.16, Theorems 2.12, 2.17, 2.21, 2.24, Propositions 2.19, 2.22,
2.23, Corollary 2.20 are similarly formulated and proved in these
cases.
\par Suppose that
$$(2.9.1)\quad 0\to {\mbox{}^1{\mathcal C}}_{\overrightarrow{u}} {\mathcal C}_{\overrightarrow{v}}
\mbox{}^2{\mathcal C}\to 0$$ is an exact sequence, where $u$ and $v$
are homomorphisms of complexes.

\par {\bf Lemma 2.10.} {\it For the exact sequence $(2.9.1)$
with a $((B,\mbox{}^2B),(B,\mbox{}^2B))$-exact homomorphism $v$ and
$\Gamma = \{ x \in {\mathcal C}: ~ dx\in Im (u) \} $ there exists a
$((\mbox{}^3B,\mbox{}^2B\otimes_{\mathcal T}
\mbox{}^1B),(\mbox{}^3B,\mbox{}^2B\otimes_{\mathcal T}
\mbox{}^1B))$-generic homomorphism $\eta : \Gamma \to
H(\mbox{}^2{\mathcal C})\times H(\mbox{}^1{\mathcal C})$ with
$\mbox{}^3B=u_{\bf \iota }(\mbox{}^1B)$ such that it is a graph of a
$((\mbox{}^2B, \mbox{}^1B),(\mbox{}^2B,\mbox{}^1 B))$-generic ${\bf
Z}$-graded homomorphism of degree $-1$ from $H(\mbox{}^2{\mathcal
C})$ to $H(\mbox{}^1{\mathcal C})$.}
\par {\bf Proof.} From $(2.9.1)$, Definition 2.8 and Remark 2.9 it follows that $u$
is a $((\mbox{}^1B,B),(\mbox{}^1B,B))$-generic homomorphism and $v$
is $((B,\mbox{}^2B),(B,\mbox{}^2B))$-generic homomorphism. From the
exactness of $(2.9.1)$, Definitions 2.1 and 2.8 it follows that
$u_{\bf \iota }: \mbox{}^1B\to B$ is injective, hence $u_{\bf \iota
}: \mbox{}^1B\to \mbox{}^3B$ is an isomorphism. Since $v_{\bf \iota
}: B\to \mbox{}^2B$ is an isomorphism, then $v_{\bf \iota }^{-1}:
\mbox{}^2B\to B$ also is an isomorphism of ${\mathcal T}$-algebras.
Therefore, $u_{\bf \iota }^{-1}\circ v_{\bf \iota }^{-1}:
\mbox{}^2B\to \mbox{}^1B$ is a homomorphism of ${\mathcal
T}$-algebras (see $(2.3.4)$). Thus there exists a homomorphism $\eta
_{\bf \iota }: \mbox{}^3B\to \mbox{}^2B\otimes_{\mathcal T}
\mbox{}^1B$.
\par Then $d(u^{-1}(dx))=u^{-1}(dd(x))=0$, consequently, $u^{-1}(dx)\in
{\mathcal Z}(\mbox{}^1{\mathcal C})$. On the other hand,
$dv(x)=v(dx)\in Im (v\circ u)=0$, hence $v(x)\in {\mathcal
Z}(\mbox{}^2{\mathcal C})$. For each $x\in \Gamma $ let $\eta (x) =
(v(x), u^{-1}(dx))$. \par Let $x\in \Gamma $ and $v(x)\in {\mathcal
B}(\mbox{}^2{\mathcal C})$. Therefore, there exists $\mbox{}^2y\in
\mbox{}^2{\mathcal C}$ such that $v(x)=d(\mbox{}^2y)$. Then there is
$y\in {\mathcal C}$ such that $\mbox{}^2y=v(y)$, consequently,
$v(x)=v(dy)$. Then there exists $\mbox{}^1b\in \mbox{}^1{\mathcal
C}$ such that $x-dy=u(\mbox{}^1b)$, hence $dx=u(d(\mbox{}^1b))$ and
$u^{-1}(dx)=d(\mbox{}^1b)$ belongs to ${\mathcal
B}(\mbox{}^1{\mathcal C})$. Thus $u^{-1}(dx)\in {\mathcal
B}(\mbox{}^1{\mathcal C})$.
\par Notice that for each element $\mbox{}^2y\in {\mathcal Z}(\mbox{}^2{\mathcal C})$ there exists $x\in {\mathcal C}$ such that
$\mbox{}^2y=v(x)$ and $v(dx)=0$. This means that $dx\in Im (u)$,
hence $x\in \Gamma $. The homomorphism $\eta $ is bihomogeneous of
degree $(0,-1)$. Thus $\eta (\Gamma )$ is a graph of the
$((\mbox{}^2B, \mbox{}^1B),(\mbox{}^2B, \mbox{}^1B))$-generic
homomorphism of degree $-1$ . This homomorphism is ${\bf Z}$-graded,
because complexes ${\mathcal C}$, $\mbox{}^1{\mathcal C}$,
$\mbox{}^2{\mathcal C}$ and homomorphisms $u$, $v$ are ${\bf
Z}$-graded.

\par {\bf Definition 2.11.} The $((\mbox{}^2B, \mbox{}^1B),(\mbox{}^2B, \mbox{}^1B))$-generic ${\bf Z}$-graded
homomorphism from Lemma 2.10 of degree $-1$ from
$H(\mbox{}^2{\mathcal C})$ to $H(\mbox{}^1{\mathcal C})$ is called a
connecting homomorphism relative to the exact sequence $(u,v)$ and
it is denoted by $\partial (u,v)$ or $\partial _{(u,v)}$, where
$\partial _n (u,v): H_n(\mbox{}^2{\mathcal C})\to
H_{n-1}(\mbox{}^1{\mathcal C})$ and $\partial ^n (u,v):
H^n(\mbox{}^2{\mathcal C})\to H^{n-1}(\mbox{}^1{\mathcal C})$ are
its homogeneous components.

\par {\bf Theorem 2.12.} {\it For the exact sequence $(2.9.1)$ of
complexes with a $((B,\mbox{}^2B),(B,\mbox{}^2B))$-exact
homomorphism $v$ there exists an exact sequence
\par $(2.12.1)$ $...\to H_{n+1}(\mbox{}^2{\mathcal C})_{\overrightarrow{\partial _{n+1}(u,v)}}H_n(\mbox{}^1{\mathcal C})
_{\overrightarrow{H_{n}(u)}}H_n({\mathcal
C})_{\overrightarrow{H_{n}(v)}}$\par $H_n(\mbox{}^2{\mathcal
C})_{\overrightarrow{\partial _{n}(u,v)}}H_{n-1}(\mbox{}^1{\mathcal
C})\to ...$
\par with a $((\mbox{}^2B, \mbox{}^1B),(\mbox{}^2B, \mbox{}^1B))$-generic homomorphism $\partial _{n}(u,v)$,
\par a $((\mbox{}^1B, B),(\mbox{}^1B, B))$-generic homomorphism $H_n(u)$ and \par a
$((B,\mbox{}^2B),(B,\mbox{}^2B))$-generic homomorphism $H_n(v)$ for
each $n$. }
\par {\bf Proof.} From Lemma 2.10 it follows that the homomorphism
\par $(2.12.2)$ $\partial _{n}(u,v): H_n(\mbox{}^2{\mathcal C})\to H_{n-1}(\mbox{}^1{\mathcal C})$ \\ exists and it is
$((\mbox{}^2A, \mbox{}^1A),(\mbox{}^2A, \mbox{}^1A))$-generic for
each $n$. Then the exact sequence $(2.9.1)$ induces the exact
sequences
\par $(2.12.3)$ $0\to {\mathcal Z}_n(\mbox{}^1{\mathcal C})_{\overrightarrow{{\mathcal Z}_{n}(u)}}{\mathcal Z}_n({\mathcal C})
_{\overrightarrow{{\mathcal Z}_{n}(v)}}{\mathcal
Z}_n(\mbox{}^2{\mathcal C})$ and
\par $(2.12.4)$ $\mbox{}^1{\mathcal C}_n/{\mathcal B}_n(\mbox{}^1{\mathcal C})_{\overrightarrow{\bar{u}_n}}
{\mathcal C}_n/{\mathcal B}_n({\mathcal C})_{\overrightarrow{\bar{v}_n}}\mbox{}^2{\mathcal C}_n/{\mathcal B}_n(\mbox{}^1{\mathcal C})\to 0$ \\
with $((\mbox{}^1B, B),(\mbox{}^1B, B))$-generic homomorphisms
${\mathcal Z}_n(u)$ and $\bar{u}_n$, $((B,\mbox{}^2B),(B,\mbox{}^2
B))$-generic homomorphisms ${\mathcal Z}_n(v)$ and $\bar{v}_n$ for
each $n$. On the other hand, the sequence
\par $(2.12.5)$ $0\to H_n({\mathcal C})_{\overrightarrow{i_n}} {\mathcal C}_n/{\mathcal B}_n({\mathcal C})
_{\overrightarrow{d_n}} {\mathcal Z}_{n-1}({\mathcal
C})_{\overrightarrow{p_{n-1}}} H_{n-1}({\mathcal C})\to 0$ \\ is
exact with $B$-exact homomorphisms $i_n$, $d_n$, $p_{n-1}$ for each
$n$, similarly for $\mbox{}^1{\mathcal C}$ and $\mbox{}^2{\mathcal
C}$. Then the exact sequences $(2.12.3)$ and $(2.12.4)$ induce the
exact sequence
\par $(2.12.6)$ $H_n(\mbox{}^1{\mathcal C})_{\overrightarrow{\hspace{0.5cm}H_n(u)\hspace{0.5cm}}} H_n({\mathcal C})
_{\overrightarrow{\hspace{0.5cm}H_n(v)\hspace{0.5cm}}}
H_n(\mbox{}^2{\mathcal C})$
\\ with a $((\mbox{}^1B, B),(\mbox{}^1B, B))$-generic homomorphism $H_n(u)$ and
a $((B, \mbox{}^2B),(B, \mbox{}^2B))$-generic homomorphism $H_n(v)$
for each $n$. Note that the homomorphism $\partial _{n}(u,v)$ is
obtained from the $((\mbox{}^2B, \mbox{}^1B),(\mbox{}^2B,
\mbox{}^1B))$-generic homomorphism $u_{n-1}^{-1}\circ d_n\circ
v_n^{-1}: {\mbox{}^2{\mathcal C}}_n\to \mbox{}^1{\mathcal C}_{n-1}$
by restricting on ${\mathcal Z}_n(\mbox{}^2{\mathcal C})$ and
${\mathcal Z}_{n-1}(\mbox{}^1{\mathcal C})$ and then using quotient
maps onto $H_n(\mbox{}^2{\mathcal C})$ and
$H_{n-1}(\mbox{}^1{\mathcal C})$ correspondingly by the construction
in Lemma 2.10. Therefore the exact sequences of the types
$(2.12.4)$-$(2.12.6)$ imply that there exists the exact sequence
$(2.12.1)$ such that the the homomorphism from
$H_n(\mbox{}^2{\mathcal C})$ into $H_{n-1}(\mbox{}^1{\mathcal C})$
coincides with $\partial _{n}(u,v)$ in $(2.12.2)$.

\par {\bf Definition 2.13.} We consider the cartesian product $X\times Y$
of $G$-graded $B$-bimodlues $X$ and $Y$ (see Definition 2.3). Let
$X\times _BY$ be a $G$-graded $B$-bimodule generated from $X\times
Y$ using finite additions of elements $(x,y)\in X\times Y$ and the
left and right multiplications on elements $a\in B$ such that
\par $(2.13.1)$ $(x,y)+(x_1,y_1)=(x+x_1,y+y_1)$ and \par $(2.13.2)$ $a(x,y)=(ax,ay)$ and
$(x,y)a=(xa,ya)$ and \par $(2.13.3)$ $g(X_e,Y_e)=(X_g,Y_g)$ and
$(X_e,Y_e)g=(X_g,Y_g)$ (see also $(2.1.1)$) for each $x$ and $x_1$
in $X$, $~y$ and $y_1$ in $Y$, $ ~ a\in B$, $~g\in G$.
\par Suppose that $X$, $Y$ and $Z$ are $G$-graded $B$-bimodules.
\par $(2.13.4)$. Let $\Lambda  : X\times Y\to Z$ be a ${\mathcal C}(B)$-bilinear
map. Let also $\Lambda $ satisfy the following identities: \par
$(2.13.5)$ $\Lambda  (x_gb_h,y_s)={\sf t}_3(g,h,s)\Lambda
(x_g,b_hy_s)$ and $\Lambda (cx,y)=c\Lambda (x,y)$ and
\par $\Lambda (x,yc) =\Lambda (x,y)c$ for each $c\in N(B)$, $x\in X$, $y\in
Y$, $g$ and $h$ and $s$ in $G$. If $\Lambda $ fulfills conditions
$(2.13.4)$ and $(2.13.5)$, then it will be said that the map
$\Lambda $ is $G$-balanced.
\par  Let $C$ be a $G$-graded $B$-bimodule supplied
with a ${\mathcal C}(B)$-bilinear map $\xi : X\times Y\to C$ denoted
by $\xi (x,y)=x\otimes y$ for each $x\in X$ and $y\in Y$ such that
\par $(2.13.6)$ $X\otimes _AY$ is generated by a set $ \{ x\otimes y:
~ x\in X, ~ y\in Y \} $ and
\par $(2.13.7)$ if $\Lambda : X\times Y\to Z$ is a $G$-balanced map
of $G$-graded $B$-bimodules $X$, $Y$ and $Z$, and for each fixed
$x\in X$ the map $\Lambda (x, \cdot ): Y\to Z$ and for each fixed
$y\in Y$ the map $\Lambda (\cdot , y): X\to Z$ are $G$-graded
homomorphisms of $G$-graded $B$-bimodules, then there exists a
$G$-graded homomorphism $\psi : C\to Z$ of $G$-graded $B$-bimodules
such that $\psi (x\otimes y)=\Lambda (x,y)$ for each $x\in X$ and
$y\in Y$.
\par If conditions $(2.13.6)$ and $(2.13.7)$ are satisfied, then the
$G$-graded $B$-bimodule $C$ is called a $G$-smashed tensor product
(or shortly tensor product) of $X$ with $Y$ over $B$ and denoted by
$X\otimes _BY$.
\par Similarly if $X$ is the $G$-graded $B$-bimodule (or right $B$-module), $Y$ is the
$G$-graded left $B$-module (or the $B$-bimodule), then the
$G$-smashed tensor product $X\otimes _BY$ of $X$ with $Y$ over $B$
is defined and it is the $G$-graded left $B$-module (or right
$B$-module correspondingly).

\par {\bf Definition 2.14.} A $G$-graded $B$-bimodule $X$ is called flat if for
each exact sequence of $G$-graded right $B$-modules $Y$,
$\mbox{}^1Y$, $\mbox{}^2Y$ and $B$-epigeneric homomorphisms $u$,
$v$:
\par $(2.14.1)$  ${\mbox{}^1Y}_{\overrightarrow{\hspace{0.5cm}u\hspace{0.5cm}}} Y
_{\overrightarrow{\hspace{0.5cm}v\hspace{0.5cm}}} \mbox{}^2Y$
\par a sequence of ${\bf Z}$-linear homomorphisms
\par $(2.14.2)$  ${\mbox{}^1Y\bigotimes _BX}_{\overrightarrow{\hspace{0.5cm}u\otimes 1\hspace{0.5cm}}}
Y\bigotimes _BX _{\overrightarrow{\hspace{0.5cm}v\otimes
1\hspace{0.5cm}}} \mbox{}^2Y\bigotimes _BX$
\par is exact, where $(u_{\bf \iota }\otimes 1)(B\bigotimes
_{\bf Z}B)=B\bigotimes _{\bf Z}B$ and $(v_{\bf \iota }\otimes
1)(B\bigotimes _{\bf Z}B)=B\bigotimes _{\bf Z}B$.

\par {\bf Proposition 2.15.} {\it A $G$-graded $B$-bimodule $X$ is flat if and
only if for each injective $B$-epigeneric homomorphism $u:
\mbox{}^1Y\to Y$ of $G$-graded right $B$-modules the ${\bf
Z}$-linear homomorphism $u\otimes 1: \mbox{}^1Y\bigotimes _BX\to
Y\bigotimes _BX$ is injective and $B$-epigeneric.}
\par {\bf Proof.} If the module $X$ is flat and a
homomorphism $u: \mbox{}^1Y\to Y$ is $B$-epigeneric and injective,
thence the following sequence
\par $(2.15.1)$  $0\rightarrow {\mbox{}^1Y}_{\overrightarrow{\hspace{0.5cm}u\hspace{0.5cm}}} Y$
\\ is exact and consequently, the homomorphism $u\otimes 1$ is
${\bf Z}$-linear and injective, where $(u_{\bf \iota }\otimes
1)(B\bigotimes _{\bf Z}B)=B\bigotimes _{\bf Z}B$.
\par If the sequence $(2.14.1)$ is exact with $B$-epigeneric
homomorphisms, then we put $Z_2=v(Y)$. Let $i: \mbox{}^2Z\to
\mbox{}^2Y$ be a canonical embedding and let $p: Y\to \mbox{}^2Z$ be
such that $v(y)$ corresponds to $y\in Y$. Then $i$ and $p$ are
$B$-epigeneric, because $v_{\bf \iota }(B)=B$. Therefore, the
following sequence
\par $(2.15.2)$  ${\mbox{}^1Y}_{\overrightarrow{\hspace{0.5cm}u\hspace{0.5cm}}} Y
_{\overrightarrow{\hspace{0.5cm}p\hspace{0.5cm}}}
{Z_2}_{\overrightarrow{\hspace{0.5cm}\xi \hspace{0.5cm}}} 0$
\\ is exact, where $\xi (Z_2)=(0)$. Hence the following sequence
\par $(2.15.3)$  ${\mbox{}^1Y\bigotimes _BX}_{\overrightarrow{\hspace{0.5cm}u\otimes 1\hspace{0.5cm}}}
Y\bigotimes _BX _{\overrightarrow{\hspace{0.5cm}p\otimes
1\hspace{0.5cm}}} Z_2\bigotimes _BX$ \\ also is exact with ${\bf
Z}$-linear homomorphisms $u\otimes 1$ and $p\otimes 1$, where
$(u_{\bf \iota }\otimes 1)(B\bigotimes _{\bf Z}B)=B\bigotimes _{\bf
Z}B$ and $(p_{\bf \iota }\otimes 1)(B\bigotimes _{\bf
Z}B)=B\bigotimes _{\bf Z}B$. Then $v=i\circ p$, consequently,
$v\otimes 1=(i\otimes 1)\circ (p\otimes 1)$ is ${\bf Z}$-linear with
$(v_{\bf \iota }\otimes 1)(B\bigotimes _{\bf Z}B)=B\bigotimes _{\bf
Z}B$. Since $i\otimes 1$ is injective, then $Ker (v\otimes 1)=Ker
(p\otimes 1)=Im (u\otimes 1)$, consequently, the sequence $(2.14.2)$
is exact.

\par {\bf Definition 2.16.} Let $({\mathcal C},d)$ be a $G$-graded $B$-complex and
$(\mbox{}^1{\mathcal C},\mbox{}^1d)$ be a $\mbox{}^1G$-graded
$\mbox{}^1B$-complex and let $f$ and $g$ be two
$((B,\mbox{}^1B),(B,\mbox{}^1B))$-generic homomorphisms of
${\mathcal C}$ into $\mbox{}^1{\mathcal C}$. A
$((B,\mbox{}^1B),(B,\mbox{}^1B))$-generic homomorphism $s$ of ${\bf
Z}$-degree $1$ from ${\mathcal C}$ into $\mbox{}^1{\mathcal C}$ such
that
\par $(2.16.1)$ $\quad g-f=\mbox{}^1d\circ s+s\circ d$ \\ is called a
homotopy relating $f$ with $g$. It is said that the homomorphisms
$f$ and $g$ are homotopic.

\par {\bf Proposition 2.17.} {\it If $f$ and $g$ are homotopic
$((B,\mbox{}^1B),(B,\mbox{}^1B))$-generic homomorphisms of
${\mathcal C}$ into $\mbox{}^1{\mathcal C}$ (see Definition 2.16),
then $H(f)=H(g)$.}
\par {\bf Proof.} If $s$ is a homotopy relating $f$ with $g$, then
$(g-f)({\mathcal Z}({\mathcal C}))=(\mbox{}^1d\circ s+s\circ
d)({\mathcal Z}({\mathcal C}))=(\mbox{}^1d\circ s)({\mathcal
Z}({\mathcal C}))\subset {\mathcal B}(\mbox{}^1{\mathcal C}) $ by
Conditions $(2.3.3)$-$(2.3.5)$, because $d$ is $B$-generic and
$\mbox{}^1d$ is $\mbox{}^1B$-generic, hence $H(g-f)=0$ and
consequently, $H(g)=H(f)$.

\par {\bf Lemma 2.18.} {\it Assume that $X$ is a $\mbox{}^1G^n$-graded left
$\mbox{}^1B$-module, $Y$ is a $\mbox{}^2G^n$-graded left
$\mbox{}^2B$-module, $f: X\to Y$ is a
$(\mbox{}^1B,\mbox{}^2B)$-generic homomorphism, $n\in {\bf N}$ (see
Definitions 2.3, 2.4). If $f_{\bf \iota }: \mbox{}^1B\to \mbox{}^2B$
is injective, then $f_{\bf \iota }(\mbox{}^1B)$ is isomorphic with
$\mbox{}^1B$. If $f_{\bf \iota }: \mbox{}^1B\to \mbox{}^2B$ is
bijective onto (i.e. injective and surjective), then $\mbox{}^1B$
and $\mbox{}^2B$ are isomorphic.}
\par {\bf Proof.} From Conditions $(2.3.3)$ and $(2.3.5)$ it follows
that $f_{\bf \iota }(\mbox{}^1G)$ is isomorphic with $\mbox{}^1G$,
if $f_{\bf \iota }: \mbox{}^1B\to \mbox{}^2B$ is injective. On the
other hand, $f$ and $f_{\bf \iota }$ are the left ${\mathcal
T}$-homomorphisms. Therefore, $f_{\bf \iota }(\mbox{}^1A)$ is
isomorphic with $\mbox{}^1A={\mathcal T}[\mbox{}^1G]$, where
$\mbox{}^1A$ is the metagroup algebra (see Definition A2). There
exists a univalent left ${\mathcal T}$-homomorphism $f_{\bf \iota
}^{-1}: f_{\bf \iota }(\mbox{}^1B)\to \mbox{}^1B$, because $f_{\bf
\iota }$ is injective. Then Conditions $(2.3.4)$ and $(2.3.5)$ imply
that $f_{\bf \iota }(\mbox{}^1B)$ is isomorphic with $\mbox{}^1B$.
\par Therefore, if $f_{\bf \iota }$ is bijective from $\mbox{}^1B$ onto
$\mbox{}^2B$, then $\mbox{}^1B$ and $\mbox{}^2B$ are isomorphic as
$\mbox{}^1G^n$-graded $\mbox{}^1A$-algebra and $\mbox{}^2G^n$-graded
$\mbox{}^2A$-algebra respectively.

\par {\bf Proposition 2.19.} {\it Assume that $\mbox{}^{k}{\mathcal C}$ are $\mbox{}^kG$-graded
$\mbox{}^kB$-complexes and $f: \mbox{}^{1}{\mathcal C}\to
\mbox{}^{2}{\mathcal C}$, $g: \mbox{}^{1}{\mathcal C}\to
\mbox{}^{2}{\mathcal C}$, $\psi : \mbox{}^{3}{\mathcal C}\to
\mbox{}^{1}{\mathcal C}$, $\eta : \mbox{}^{2}{\mathcal C}\to
\mbox{}^{4}{\mathcal C}$ are
$((\mbox{}^jB,\mbox{}^kB),(\mbox{}^jB,\mbox{}^kB))$-generic
homomorphisms of complexes with $(j,k)=(1,2)$ for $f$ and $g$,
$(j,k)=(3,1)$ for $\psi $, $(j,k)=(2,4)$ for $\eta $. If $s$ is a
homotopy relating $f$ with $g$, then $\eta \circ s\circ \psi $ is a
homotopy relating $\eta \circ f\circ \psi $ with $\eta \circ g\circ
\psi $.}
\par {\bf Proof.} The composition $\eta \circ s\circ \psi $ is
a $((\mbox{}^3B,\mbox{}^4B),(\mbox{}^3B,\mbox{}^4B))$-generic
homomorphisms of ${\bf Z}$-degree $1$ from $\mbox{}^3{\mathcal C}$
into $\mbox{}^4{\mathcal C}$. Then using Definitions 2.8 and 2.16
one verifies the assertion of this proposition.

\par {\bf Corollary 2.20.} {\it Let $\mbox{}^{k}{\mathcal C}$ be $\mbox{}^kG$-graded
$\mbox{}^kB$-complexes, where $k\in \{ 1, 2, 3 \} $, let also $f^j$
and $g^j$ be
$((\mbox{}^jB,\mbox{}^{j+1}B),(\mbox{}^jB,\mbox{}^{j+1}B))$-generic
homomorphisms from  $\mbox{}^{j}{\mathcal C}$ into
$\mbox{}^{j+1}{\mathcal C}$ for $j=1$ and $j=2$. Let $s^j$ be a
homotopy of $f^j$ with $g^j$ for $j=1$ and $j=2$. Then $s^2\circ
f^1+ g^2\circ s^1$ is a homotopy relating $f^2\circ f^1$ with
$g^2\circ g^1$.}
\par {\bf Proof.} From the conditions of this Corollary and Definition 2.3
it follows that the homomorphisms $s^2\circ f^1$, $f^2\circ f^1$,
$g^2\circ f^1$, $g^2\circ s^1$, $g^2\circ f^1$ and $g^2\circ g^1$
are $((\mbox{}^1B,\mbox{}^{3}B),(\mbox{}^1B,\mbox{}^{3}B))$-generic.
In view of Proposition 2.19 $s^2\circ f^1$ relates $f^2\circ f^1$
with $g^2\circ f^1$, while $g^2\circ s^1$ relates $g^2\circ f^1$
with $g^2\circ g^1$, hence $s^2\circ f^1+ g^2\circ s^1$ relates
$f^2\circ f^1$ with $g^2\circ g^1$.

\par {\bf Definition 2.21.} A $((\mbox{}^1B,\mbox{}^{2}B),(\mbox{}^1B,\mbox{}^{2}B))$-generic
homomorphism $f^1: \mbox{}^{1}{\mathcal C}\to \mbox{}^{2}{\mathcal
C}$ of $\mbox{}^kG$-graded $\mbox{}^kB$-complexes, where $k\in \{ 1,
2 \} $, is called a homotopism if there exists a
$((\mbox{}^2B,\mbox{}^{1}B),(\mbox{}^2B,\mbox{}^{1}B))$-generic
homomorphism $f^2: \mbox{}^{2}{\mathcal C}\to \mbox{}^{1}{\mathcal
C}$ such that $f^2\circ f^1$ and $f^1\circ f^2$ are homotopic to
$1_{\mbox{}^1{\mathcal C}}$ and $1_{\mbox{}^2{\mathcal C}}$
respectively. The complex $\mbox{}^1{\mathcal C}$ is homotopic to
$0$, if $1_{\mbox{}^1{\mathcal C}}$ is homotopic to
$0_{\mbox{}^1{\mathcal C}}$.

\par {\bf Proposition 2.22.} {\it If $f^1$ is a homotopism, then it is a
homologism. Moreover, if $g^1$ is homotopic to $f^1$, then $g^1$
also is a homotopism.}
\par {\bf Proof.} In the notation of Definition 2.21
$H(f^2)\circ H(f^1)=H(f^2\circ f^1)=H(1_{\mbox{}^1{\mathcal
C}})=1_{H(\mbox{}^1{\mathcal C})}$ by Proposition 2.17. Similarly
$H(f^1)\circ H(f^2)=1_{H(\mbox{}^2{\mathcal C})}$. Hence $H(f^1)$ is
bijective and $f^1$ is a homologism (see Remark 2.9).
\par Then for $f^1$ and $g^1$ one gets that $(f^2\circ g^1)$ is homotopic to
$(f^2\circ f^1)$, consequently, to $1_{\mbox{}^1{\mathcal C}}$.
Analogously $(g^1\circ f^2)$ is homotopic to $(f^1\circ f^2)$,
consequently, to $1_{\mbox{}^2{\mathcal C}}$ by Proposition 2.19.
Thus $g^1$ is the homotopism.

\par {\bf Proposition 2.23.} {\it Suppose that $({\mathcal C},d)$ is a $G$-graded
$B$-complex (see Definition 2.8). Then the following conditions are
equivalent:
\par $(2.23.1)$ there exists a $B$-generic homotopism of $({\mathcal
C},d)$ onto $(H({\mathcal C}),0)$;
\par $(2.23.2)$ there exists a $B$-generic and ${\bf Z}$-graded of degree $1$
endomorphism $s$ of the module ${\mathcal C}$ for which $d=d\circ
s\circ d$;
\par $(2.23.3)$ ${\mathcal B}({\mathcal C})$ and ${\mathcal Z}({\mathcal C})$ are direct multipliers of ${\mathcal C}$;
\par $(2.23.4)$ $({\mathcal C},d)$ is a direct sum of subcomplexes, which have length either $0$
or $1$ and zero homology.}
\par {\bf Proof.} $(2.23.1)$ $\Rightarrow $ $(2.23.2)$. Assume that
$f: {\mathcal C}\to H({\mathcal C})$ is a homotopism. Therefore,
there exists a morphism of complexes $g: H({\mathcal C})\to
{\mathcal C}$ and the $B$-generic and ${\bf Z}$-graded of degree $1$
endomorphism $s$ of the module ${\mathcal C}$ such that $g\circ
f=1_{\mathcal C}-s\circ d-d\circ s$. From $d\circ g=0$ and $g\circ
0=0$ it follows that $d\circ g\circ f=0$ and consequently, $d-d\circ
s\circ d-d\circ d\circ s=d-d\circ s\circ d=0$. The latter implies
$(2.23.2)$.
\par $(2.23.2)$ $\Rightarrow $ $(2.23.3)$. Let $s$ be an
endomorphism provided by Condition $(2.23.2)$. Hence $d\circ
(1_{\mathcal C}-s\circ d)=0$. The maps $d$ and $s$ are $B$-generic,
hence $1_{\mathcal C}-s\circ d$ is the ${\mathcal T}$-linear
$B$-generic projector from ${\mathcal C}$ onto ${\mathcal
Z}({\mathcal C})$. On the other hand, $d\circ s\circ d=d$,
consequently, $d\circ s$ is the ${\mathcal T}$-linear $B$-generic
projector from ${\mathcal C}$ onto ${\mathcal B}({\mathcal C})$.
\par  $(2.23.3)$ $\Rightarrow $ $(2.23.4)$. For each integer $n$ we
consider ${\mathcal Z}_n={\mathcal Z}_n({\mathcal C})$ and
${\mathcal B}_n={\mathcal B}_n({\mathcal C})$. Choose
$G^{|n|+2}$-graded $B$-sub-bimodules ${\mathcal P}_n$ and ${\mathcal
Q}_n$ in ${\mathcal C}_n$ for which ${\mathcal C}_n={\mathcal
Z}_n\oplus {\mathcal Q}_n$ and ${\mathcal Z}_n={\mathcal B}_n\oplus
{\mathcal P}_n$. By $({\mathcal C}(p),d(p))$ we denote the $p$-th
translate of $({\mathcal C},d)$, where ${\mathcal C}(p)_n={\mathcal
C}_{n+p}$ and ${\mathcal C}(p)^n={\mathcal C}^{n-p}$,
$d(p)=(-1)^pd$. We take ${\mathcal S}_{(n)}:={\mathcal P}_n(-n)$ and
${\mathcal T}_{(n)}={\mathcal Q}_n(-n)\oplus {\mathcal
B}_{n-1}(1-n)$ subcomplexes in $({\mathcal C},d)$, because
${\mathcal Z}({\mathcal C}(p))={\mathcal Z}({\mathcal C})(p)$,
${\mathcal B}({\mathcal C}(p))={\mathcal B}({\mathcal C})(p)$ and
$H({\mathcal C}(p))=H({\mathcal C})(p)$. This gives $({\mathcal
C},d)=\bigoplus _{n\in {\bf Z}} ({\mathcal S}_{(n)}\oplus {\mathcal
T}_{(n)})$. Note that for each $n\in {\bf Z}$ the complex ${\mathcal
S}_{(n)}$ is either nil or of null length, while ${\mathcal
T}_{(n)}$ is either nil or of length one with zero homology. The
latter implies $(2.23.4)$. \par $(2.23.4)$ $\Rightarrow $
$(2.23.1)$. Condition $(2.23.1)$ is satisfied if the complex
${\mathcal C}$ is of null length or of length one with zero
homology.
\par {\bf Definition 2.24.} It is said that a $G$-graded $B$-complex
$({\mathcal C},d)$ is split, if it satisfies equivalent conditions
of Proposition 2.23. A $B$-generic endomorphism $s$ of ${\mathcal
C}$ satisfying $(2.23.2)$ is called a splitting of ${\mathcal C}$.
\par Let $X$ be a $G$-graded $B$-bimodule, and let $({\mathcal P},d)$ be
a $G$-graded $B$-complex and let $({\mathcal P},d)$ be null from the
right (or left) and let $p: {\mathcal P}\to X$ be a homologism. Then
the pair $({\mathcal P},p)$ (or $(p, {\mathcal P})$ respectively) is
called a left (or right respectively) $G$-graded resolution of $X$.
A length of the $G$-graded $B$-complex $({\mathcal P},d)$ is called
a length of the resolution. If $({\mathcal P},p)$ and
$(\mbox{}^1{\mathcal P},\mbox{}^1p)$ are two left resolutions (or
two right resolutions $(p, {\mathcal P})$ and $(\mbox{}^1p,
\mbox{}^1{\mathcal P})$) and $f: {\mathcal P}\to \mbox{}^1{\mathcal
P}$ is a $B$-generic or $B$-exact morphism of complexes such that
$\mbox{}^1p\circ f=p$ (or $f\circ p=\mbox{}^1p$ respectively), then
$f$ is called a $B$-generic or $B$-exact respectively morphism of
resolutions.
\par Analogously are considered complexes and splittings in
other cases: if ${\mathcal C}$ is a $G$-graded left $B$-module or
right $B$-module and $X$ is a $G$-graded left $B$-module or right
$B$-module correspondingly.

\par {\bf Remark 2.25.} Let $A={\mathcal T}[G]$ be a metagroup algebra over a commutative
associative unital ring ${\mathcal T}$, $B$ be a $G$-graded
$A$-algebra, where $G$ is a (nonassociative) metagroup. We put \par
$(2.25.1)$ ${\mathcal K}_n=0$ for each $n<-1$, \par $(2.25.2)$
${\mathcal K}_{-1}=B$, ${\mathcal K}_0 = B \otimes_{\mathcal T}B$
and by induction \par $(2.25.3)$ ${\mathcal K}_{n+1} = {\mathcal
K}_n\otimes_{\mathcal T}B$ for each natural number $n$ (see
Definition 2.13). Therefore ${\mathcal K}_n$ is supplied with a
two-sided $B$-module structure.

\par {\bf Proposition 2.26.} {\it The $B$-bimodule ${\mathcal K}_n$ is
$G^{|n|+2}$-graded for each $n\in {\bf Z}$. Moreover, if $G$,
${\mathcal T}$ and $B$ are nontrivial, then ${\mathcal K}_n$ is
essentially $G^{n+2}$-graded for each $n\ge -1$.}
\par {\bf Proof.} If $n<0$ and ${\mathcal K}_n=0$, then it can be supplied with the
trivial $G^{|n|+2}$-gradation such that $({\mathcal K}_n)_{ \{
g_0,...,g_{|n|+1} \} _{l(n+2)}}=0$ for each $g_0,...,g_{|n|+1}$ in
$G$. If $n<0$ and ${\mathcal K}_n=B$, then it is $G$-graded and by
Lemma 2.7 it can be supplied with the $G^{|n|+2}$-graded structure.
\par For each $n\in {\bf N}$ it satisfies the following identities:
\par $(2.26.1)$ $\forall p\in {\mathcal T}$, $~p\cdot (z_{g_0},...,z_{g_{n+1}})= ((pz_{g_0}),...,z_{g_{n+1}})$ and
\par $(z_{g_0},...,(z_{g_{n+1}}p))= (z_{g_0},...,z_{g_{n+1}})\cdot p$ and \par $\forall
j\in \{ 1,...,n \}$,  $~p\cdot (z_{g_0},...,z_{g_{n+1}})=
(z_{g_0},...,(pz_{g_j}),...,z_{g_{n+1}})$ and
\par $(z_{g_0},...,(z_{g_j}p),...,z_{g_{n+1}})=(z_{g_0},...,z_{g_{n+1}})\cdot p$, \par
where $0\cdot (z_{g_1},...,z_{g_n})=0$;
\par $(2.26.2)$ $(z_gz_y)\cdot (z_{g_0},...,z_{g_{n+1}})={\sf t}_3\cdot (z_g\cdot (z_y\cdot (z_{g_0},...,z_{g_{n+1}})))$ \par with
${\sf t}_3={\sf t}_3(g,y,b)$;
\par $(2.26.3)$ ${\sf t}_3\cdot ((z_{g_0},...,z_{g_{n+1}})\cdot (z_gz_y)) = ((z_{g_0},...,z_{g_{n+1}})\cdot z_g)\cdot z_y$
with ${\sf t}_3={\sf t}_3(b,g,y)$;
\par $(2.26.4)$ $(z_g\cdot (z_{g_0},...,z_{g_{n+1}}))\cdot z_y = {\sf t}_3 \cdot (z_g\cdot ((z_{g_0},...,z_{g_{n+1}})\cdot z_y))$ with
${\sf t}_3={\sf t}_3(g,b,y)$
\par $(2.26.5)$ $z_g\cdot (z_{g_0},...,z_{g_{n+1}})=$\par $ t_{n+3} (g,g_0,...,g_{n+1};v_0(n+3);l(n+3))\cdot ((z_gz_{g_0}),z_{g_1},...,z_{g_{n+1}})$
\\ where $\{ g, g_0,...,g_{n+1} \} _{v_0(n+3)} = g \{ g_0,...,g_{n+1} \} _{l(n+2)}$,
\par $\{ g_0,...,g_{n+1} \}_{l(n+2)} = \{ g_0,...,g_n \} _{l(n+1)} g_{n+1}$, \par $\{ g_0 \}_{l(1)}= g_0$,
$ \{ g_0 g_1 \} _{l(2)} =g_0g_1$; \\ where $b= \{ g_0,...,g_{n+1} \}
_{l(n+2)}$,
\par $t_{n} (g_1,...,g_{n};u(n),w(n)) := t_{n} (g_1,...,g_{n};u(n), w(n)|id)$ \\ using shortened notation;
\par $(2.26.6)$ $(z_{g_0},...,z_{g_{n+1}})\cdot z_g= t_{n+3}(g_0,...,g_{n+1},g;l(n+3),v_{n+2}(n+3))\cdot (z_{g_0},...,z_{g_n}, (z_{g_{n+1}}z_g))$
\\ for every $g, y, g_0,...,g_{n+1}$ in $G$, $z_g\in B_g$, $z_y\in
B_y$, $z_{g_j}\in B_{g_j}$ for each $j\in \{ 0,...,n+1 \} $, where
$(z_{g_0},...,z_{g_{n+1}})$ is a shortened notation of the left
ordered tensor product
\par $(...((z_{g_0}\otimes z_{g_1})\otimes z_{g_2})...\otimes z_{g_n})\otimes z_{g_{n+1}}$,
\par $\{ g_0,...,g_{n+1},g \} _{v_{n+2}(n+3)} = \{ g_0,...,g_n,g_{n+1}g \} _{l(n+2)}$.
\par If $n$ is nonnegative, $0\le n\in {\bf Z}$, then from $(2.13.6)$
and $(2.13.7)$ it follows that
\par \par $(2.26.7)$ ${\mathcal K}_n=\sum_{g_0\in G,...,g_{n+1}\in G}({\mathcal
K}_n)_{\{ g_0,...,g_{n+1} \} _{l(n+2)}} $ \par with $({\mathcal
K}_n)_{\{ g_0,...,g_{n+1} \} _{l(n+2)}}$ consisting of all elements
$z_{\{ g_0,...,g_{n+1} \} _{l(n+2)}}$ which are sums of elements of
the form $(z_{g_0},...,z_{g_{n+1}})$ with $z_{g_j}\in B_{g_j}$ for
each $j\in \{ 0,...,n+1 \} $ (see Remark 2.25). From Identities
$(2.26.1)$-$(2.26.6)$ it follows that the $B$-bimodule ${\mathcal
K}_n$ satisfies Conditions $(2.4.1)$-$(2.4.5)$.
\par If $G$, ${\mathcal T}$ and $B$ are nontrivial, then Identities
$(2.25.2)$, $(2.25.3)$, $(2.26.1)$-$(2.26.7)$ and Definitions 2.4,
2.13 imply that ${\mathcal K}_n$ is essentially $G^{n+2}$-graded for
each $n\ge -1$.

\par {\bf Proposition 2.27.} {\it Let the algebra $B$ and the $B$-bimodules
${\mathcal K}_n$ be as in Remark 2.25. Then an acyclic left
$B$-complex $({\mathcal K}(B),\partial )$ exists.}
\par {\bf Proof.} We take the $B$-bimodules ${\mathcal K}_n$ for each $n\ge 0$
as in Remark 2.25 and ${\mathcal K}_{-1}=B$. In view of Proposition
2.26 the $B$-bimodule ${\mathcal K}_n$ is $G^{|n|+2}$-graded for
each $n\in {\bf Z}$. The metagroup algebra $A$ is unital and the
$G$-graded $A$-algebra $B$ is unital such that $A$ has the natural
embedding into $B$ as $A1_B$, where $1_B=1$ is the unit element in
$B$. Therefore ${\mathcal T}1_B\subset {\sf C}(B)\subset {\sf
C}(A)1_B$ (see Definition 2.1).
\par Then we describe a boundary ${\mathcal T}$-linear operator
$\partial_n: {\mathcal K}_n\to {\mathcal K}_{n-1}$ on ${\mathcal
K}_n$ for each natural number $n$. Using the decomposition
$(2.26.7)$ it is sufficient to give it at first on
$(z_{g_0},...,z_{g_{n+1}})$ for every $g, y, g_0,...,g_{n+1}$ in
$G$, $z_g\in B_g$, $z_y\in B_y$, $z_{g_j}\in B_{g_j}$ for each $j\in
\{ 0,...,n+1 \} $. Then it has by the ${\mathcal T}$-linearity an
extension on ${\mathcal K}_n$. Therefore we put:
\par $(2.27.1)$ $\partial_n((z_g\cdot (z_{g_0},z_{g_1},...,z_{g_n},z_{g_{n+1}}))\cdot z_y) =$\par $\sum_{j=0}^{n}
(-1)^{j}\cdot t_{n+4}(g,g_0,...,g_{n+1},y;l(n+4),u_{j+1}(n+4))$\par
$\cdot ((z_g\cdot(<z_{g_0}, z_{g_1},..., z_{g_{n+1}}
>_{j+1,n+2}))\cdot z_y) $, where
\par $(2.27.2)$ $ <z_{g_0},...,z_{g_{n+1}} >_{1,n+2} := ((z_{g_0}z_{g_1}),z_{g_2},...,z_{g_{n+1}})$,
\par $(2.27.3)$ $ <z_{g_0},...,z_{g_{n+1}} >_{2,n+2} :=  (
z_{g_0},(z_{g_1}z_{g_2}),z_{g_3},...,z_{g_{n+1}} )$,...,
\par $(2.27.4)$ $ < z_{g_0},...,z_{g_{n+1}} >_{n+1,n+2} :=  (
z_{g_0},...,z_{g_{n-1}},(z_{g_n}z_{g_{n+1}}) )$,
\par $(2.27.5)$ $\partial_0 (z_g\cdot (z_{g_0},z_{g_1}))\cdot z_y =(z_g\cdot (z_{g_0}z_{g_1}))\cdot z_y$,
\par $(2.27.6)$ $\{ g_0,g_1,...,g_{n+1} \}_{l(n+2)}
:= (...((g_0g_1)g_2)...)g_{n+1} $;
\par $(2.27.7)$ $ \{ g, g_0, ...,g_{n+1}, y \}_{u_1(n+4)} :=  (g \{ (g_0g_1),
g_2,...,g_{n+1} \}_{l(n+1)})y$,...,
\par $(2.27.8)$ $ \{ g, g_0,...,g_{n+1},y \}_{u_{n+1}(n+4)} :=  (g \{
g_0,g_1,...,(g_ng_{n+1}) \}_{l(n+1)})y$
\\ for each $g, g_0,...,g_{n+1}, y$ in $G$.
On the other hand, from formulas $(1)$ and $(2)$ in Definition 1 it
follows that $t_{n+4}(g,g_0,...,g_{n+1},y;l(n+4),u_{j+1}(n+4))
=$\par $ t_{n+2} (g_0,...,g_{n+1};l(n+2),v_{j+1}(n+2))$  for each
$j=0,...,n$, where
\par $(2.27.9)$ $ \{ g_0, ...,g_{n+1} \}_{v_1(n+2)} :=  \{ (g_0g_1),
g_2,...,g_{n+1} \}_{l(n+1)}$,...,
\par $(2.27.10)$ $ \{ g_0,...,g_{n+1} \}_{v_{n+1}(n+2)} :=  \{
g_0,g_1,...,(g_ng_{n+1}) \}_{l(n+1)}$
\\ for every $g_0,...,g_{n+1}$ in $G$.
This means that $\partial _n$ is a left and right $B$-homomorphism
of $B$-bimodules, consequently, $\partial _n$ is $B$-exact.
Particularly,
\par $(2.27.11)$ $\partial _1((z_g\cdot
(z_{g_0},z_{g_1},z_{g_2}))\cdot z_y) = (z_g\cdot ((z_{g_0}
z_{g_1}),z_{g_2}))\cdot z_y - {\sf t}_{3}(g_0,g_1,g_2) \cdot
(z_g\cdot (z_{g_0}, (z_{g_1}z_{g_2})))\cdot z_y $,
\par $(2.27.12)$ $\partial _2((z_g\cdot (z_{g_0},z_{g_1},z_{g_2},z_{g_3}))\cdot z_y)=
(z_g\cdot ((z_{g_0}z_{g_1}),z_{g_2},z_{g_3}))\cdot z_y
-t_4(g_0,...,g_3; l(4),v_2(4))\cdot ((z_g\cdot
(z_{g_0},(z_{g_1}z_{g_2}),z_{g_3}))\cdot z_y)$\par $ +
t_4(g_0,...,g_3; l(4),v_3(4)) ((z_g\cdot
(z_{g_0},z_{g_1},(z_{g_2}z_{g_3})))\cdot z_y)$.
\par Then we define a ${\mathcal T}$-linear homomorphism
${\bf s}_n: {\mathcal K}_n \to {\mathcal K}_{n+1}$, which has the
form:
\par $(2.27.13)$ ${\bf s}_n(z_{g_0},...,z_{g_{n+1}})= (1,z_{g_0},....,z_{g_{n+1}})$
\\ for every $g_0,...,g_{n+1}$ in $G$.
From formulas $(1)$ in Lemma 1 in \cite{ludkctnaax19} and the
identities $(12)$-$(14)$ in Proposition 1 in \cite{ludkctnaax19} and
$(2.1.6)$ in Definition 2.1 it follows that
\par $(2.27.14)$ ${\bf s}_n((z_{g_0},...,z_{g_{n+1}})\cdot z_y) = ({\bf s}_n(z_{g_0},...,z_{g_{n+1}}))\cdot z_y$
\\ for every $g_0,...,g_{n+1}, y$ in $G$, $z_{g_j}\in B_{g_j}$ for each $j\in \{ 0,...,n+1 \} $,
$z_y\in B_y$.
\par We put ${\bf p}_n: {\mathcal K}_{n+1}\to {\mathcal K}_n$ to be a ${\mathcal T}$-linear mapping such that
\par $(2.27.15)$ ${\bf p}_n(a\otimes b)=a\cdot b$ and ${\bf p}_n(b\otimes a)=b\cdot a$ for each
$a\in {\mathcal K}_n$ and $b\in B$. Hence formulas $(13)$ and $(14)$
in Proposition 1 in \cite{ludkctnaax19} and $(2.27.14)$, $(2.27.15)$
imply that ${\bf p}_n{\bf s}_n=id$ is the identity on ${\mathcal
K}_n$, consequently, ${\bf s}_n$ is a monomorphism.
\par Therefore, from formulas $(2.27.1)$ and $(2.27.11)$, $((2.27.14)$ it follows that
\par $(2.27.16)$ $(\partial _{n+1} {\bf s}_n + {\bf s}_{n-1} \partial _n)(z_{g_0},...,z_{g_{n+1}}) =$
\par $= \partial _{n+1} (1,z_{g_0},...,z_{g_{n+1}}) +$ \par ${\bf s}_{n-1}(\sum_{j=0}^{n}
(-1)^{j} t_{n+2}(g_0,...,g_{n+1};l(n+2),v_{j+1}(n+2))\cdot $\par $
<z_{g_0}, ..., z_{g_{n+1}} >_{j+1,n+2}))=$
\par $=\sum_{j=0}^{n+1}
(-1)^{j}t_{n+3}(1,g_0,...,g_{n+1};l(n+3),v_{j+1}(n+3))\cdot $\par $
<1,z_{g_0}, z_{g_1},..., z_{g_{n+1}} >_{j+1,n+3}+$
\par
$+ \sum_{j=0}^{n} (-1)^{j}
t_{n+2}(g_0,...,g_{n+1};l(n+2),v_{j+1}(n+2))\cdot $\par $
<1,z_{g_0}, ...,z_{g_{n+1}} >_{j+2,n+3} = (z_{g_0},...,z_{g_{n+1}})$,\\
for every $g_0$,...,$g_{n+1}$ in $G$, $z_{g_j}\in B_{g_j}$ for each
$j\in \{ 0,...,n+1 \} $.
\par Then Formulas $(2.27.14)$-$(2.27.16)$ imply the homotopy conditions
\par $(2.27.17)$ $\partial _{n+1}{\bf s}_n+{\bf s}_{n-1}\partial _n =I$ for each $n\ge 0$, \\
where $I$ denotes the identity operator on ${\mathcal K}_n$. This
leads to identities: \par $\partial _n\partial _{n+1} {\bf s}_n =
\partial _n (I-{\bf s}_{n-1}\partial _n) = \partial _n - (\partial _n {\bf s}_{n-1})\partial _n = \partial _n -
(I-{\bf s}_{n-2}\partial _{n-1})\partial _n$ and hence gives the
recurrence relation
\par $(2.27.18)$ $\partial _n\partial _{n+1} {\bf s}_n = {\bf s}_{n-2}
\partial _{n-1} \partial _n$. \par Notice that Formula $(2.27.13)$ implies that
${\mathcal K}_{n+1}$ as the left $B$-module is generated by ${\bf
s}_n{\mathcal K}_n$. Utilizing the recurrence relation $(2.27.18)$
by induction in $n$ we infer that $\partial _n \partial _{n+1}=0$
for each $n\ge 0$, since $\partial _0 \partial _1=0$ according to
formulas $(2.27.1)$ and $(2.27.5)$.
\par Let $B^e := B\otimes
_AB^{op}$ be the enveloping algebra of $B$, where $B^{op}$ denotes
an opposite algebra. The latter as an ${\sf C}(B)$-linear space is
the same, but with the multiplication $x\circ y =yx$ for each $x, y
\in B^{op}$. This permits to consider the $G^2$-graded $B$-bimodule
${\mathcal K}_0$ as $B^e$. Therefore, the mapping $\partial _0:
{\mathcal K}_0\to {\mathcal K}_{-1}$ provides the augmentation
$\epsilon : B^e\to B$.
\par Thus according to identities $(2.27.17)$ the left complex ${\mathcal K}(B)$
is acyclic:
\par $(2.27.19)$ $0\leftarrow B \mbox{ }_{
\overleftarrow{\partial _0}} {\mathcal K}_0 \mbox{
}_{\overleftarrow{\partial _1}} {\mathcal K}_1 \mbox{
}_{\overleftarrow{\partial _2}} {\mathcal K}_2 \overleftarrow ...
\mbox{ }_{\overleftarrow{\partial _n}} {\mathcal K}_n \mbox{
}_{\overleftarrow{\partial _{n+1}}} {\mathcal K}_{n+1} \leftarrow
...$.

\par {\bf Remark 2.28.}  For a $G$-graded left $B$-module $X$ let $({\mathcal K}(B,X),
d)$ be a complex such that ${\mathcal K}_n(B,X)={\mathcal
K}_n\bigotimes _BX$ for each $0\le n\in {\bf Z}$ with $d_n=\partial
_n\otimes I_X$ for each $n\ge 1$ (see Proposition 2.27 and Remark
2.25), while ${\mathcal K}_n(B,X)=0$ and $d_{n+1}=0$ for each $n<0$,
where $I_X$ denotes a unit operator on $X$, such that $I_Xx=x$ for
each $x\in X$. Let a map $\epsilon _X: {\mathcal K}_0(B,X)\to X$ be
defined by the following formula:
\par $(2.28.1)$ $\quad \epsilon _X((a\otimes b)\otimes x)=(ab)x$
\\ for each $a$ and $b$ in $B$ and $x\in X$.
Formula $(2.28.1)$ above, conditions $(A3.1)$-$(A3.3)$ in Definition
A3 and the identities $(2.27.5)$, $(2.27.11)$ imply that
\par $(2.28.2)$ $\quad \epsilon _X\circ d_1=0$.
\\ This procedure induces a ${\bf Z}$-graded homomorphism
\par $(2.28.3)$ $\quad \hat{ \epsilon }_X: {\mathcal K}(B,X)\to X$.

\par {\bf Proposition 2.29.} {\it The map
$\hat{ \epsilon }_X$ (see Remark 2.25) is a homotopism of complexes
of left $B$-modules. The complex $({\mathcal K}(B,X), d)$ splits as
a $G$-graded left $B$-complex and $({\mathcal K}(B,X), \hat{
\epsilon }_X)$ is a left resolution of the $G$-graded left
$B$-module $X$.}
\par {\bf Proof.} For each $n\ge 0$ there exists a ${\mathcal T}$-linear
map \par $(2.29.1)$ $\quad v_n \{ (b_0\otimes  ... \otimes b_{n+1} )
\} _{l(n+2)}=\{ (1\otimes b_0\otimes ... \otimes b_{n+1} ) \}
_{l(n+3)}$ \\ for each $b_0$,...,$b_{n+1}$ in $B$, where $\{
(b_0\otimes b_1 ) \} _{l(2)}= b_0\otimes b_1$ and by induction $ \{
(b_0\otimes ... \otimes b_{n+1} ) \} _{l(n+2)}= \{ (b_0\otimes ...
\otimes b_n ) \} _{l(n+1)}\otimes b_{n+1}$ for each $n\ge 2$. This
$v_n$ is a homomorphism from ${\mathcal K}_n$ into ${\mathcal
K}_{n+1}$ as right $A$-modules. Proposition 2.27, Formula
$(2.29.1)$, Lemma 2.2, Remark 2.25 and Definitions 2.1, A1, A2 imply
that
\par $(2.29.2)$ $\quad d_{n+1}\circ v_n+v_{n-1}\circ d_n=1_{{\mathcal K}_n}$
\\ for each $n\ge 1$, because $t_3(e,g_1,g_2)=e$ for each $g_1$ and $g_2$ in
$G$, where $e$ is the unit element in $G$. In particular,
\par $(2.29.3)$ $\quad d_1\circ v_0(b_0\otimes b_1)=b_0\otimes b_1-1\otimes (b_0b_1)$
\\ for each $b_0$ and $b_1$ in $B$. A map $\xi : B\to B\bigotimes
_{\mathcal T}B$ such that $\xi (b)=1\otimes b$ for each $b\in B$
induces a ${\mathcal T}$-linear homomorphism $\hat{\xi }: B\to
{\mathcal K}(B)$. This implies that $\hat{\epsilon }_B\circ \hat{\xi
}=I_B$. From Formulas $(2.29.2)$ and $(2.29.3)$ it follows that
$d\circ v+v\circ d=I_{{\mathcal K}(B)}-\hat {\xi }\circ \hat
{\epsilon }_B$. Then defining $\hat {\xi }_X=\hat{\xi }\otimes I_X$,
$d_X=d\otimes I_X$, $v_X=v\otimes I_X$ we infer that $\hat{\epsilon
}_X\circ \hat{\xi }_X=I_X$ and $d_X\circ v_X+v_X\circ
d_X=I_{{\mathcal K}(B,X)}-\hat {\xi }_X\circ \hat {\eta }_X$. The
homomorphisms $d_X$, $v_X$, $\hat {\epsilon }_X$ and $\hat {\xi }_X$
are $B$-generic, since the homomorphisms $d$, $v$, $\hat {\epsilon
}$ and $\hat {\xi }$ are $B$-generic. Thus $\hat {\epsilon }_X$ is a
homotopism (see Definition 2.21). Proposition 2.27 implies that
$({\mathcal K}(B,X), \hat{ \epsilon }_X)$ is a left resolution of
the $G$-graded left $B$-module $X$, because $\hat {\epsilon }_X$ is
the homotopism.

\par {\bf Definition 2.30.} The left resolvent $({\mathcal K}(B,X),\hat
{\epsilon }_X)$ for $X$ is called the standard resolvent of the
$G$-graded left $B$-module $X$.

\section{Smashed torsion product.}
\par {\bf Definition 3.1.} Let $G$ be a metagroup, ${\mathcal T}$ be a
commutative associative ring, $A={\mathcal T}[G]$ be a metagroup
algebra, $B$ be a $G$-graded unital $A$-algebra. Let also
$({\mathcal C},d)$ and $(\mbox{}^1{\mathcal C},\mbox{}^1d)$ be
$G$-graded $B$-complexes, where either ${\mathcal C}$ is the
$B$-bimodule and  $\mbox{}^1{\mathcal C}$ is the left $B$-module or
${\mathcal C}$ is the right $B$-module and  $\mbox{}^1{\mathcal C}$
is the $B$-bimodule (see also Definitions 2.8, 2.13 and Remark 2.9).
Let a $G$-smashed tensor product ${\mathcal C}\bigotimes _
B\mbox{}^1{\mathcal C}$ be supplied also with the ${\bf
Z}$-gradation
\par $(3.1.1)$ $\quad ({\mathcal C}\bigotimes  _ B\mbox{}^1{\mathcal C})_n=\sum_{j+l=n} ({\mathcal C}_j\bigotimes _
B\mbox{}^1{\mathcal C}_l)$.
\par We put $D$ to be a ${\mathcal T}$-linear endomorphism of degree
$-1$ on ${\mathcal C}\bigotimes _B\mbox{}^1{\mathcal C}$ such that
\par $(3.1.2)$ $\quad D(x\otimes \mbox{}^1x)=(dx)\otimes \mbox{}^1x+(-1)^jx\otimes (\mbox{}^1dx)$
\\ for each $x\in {\mathcal C}_j$ and $\mbox{}^1x\in \mbox{}^1{\mathcal C}_i$,
$i $ and $j$ in ${\bf Z}$. The $G$-graded $B$-complex $({\mathcal
C}\bigotimes _ B\mbox{}^1{\mathcal C}, D)$ is called a $G$-smashed
tensor product (or shortly tensor product) of the complexes
$({\mathcal C},d)$ and $(\mbox{}^1{\mathcal C},\mbox{}^1d)$.

\par {\bf Remark 3.2.} In view of Definitions 2.8, 3.1 and A3
\par $(3.2.1)$ $\quad D\circ D=0$, \\ since
$\quad D\circ D(x\otimes \mbox{}^1x)=(d\circ dx)\otimes \mbox{}^1x
+(-1)^{j-1}(dx)\otimes (\mbox{}^1d\mbox{}^1x)$\par $ +
(-1)^j(dx)\otimes
(\mbox{}^1d\mbox{}^1x) + x\otimes (\mbox{}^1d\circ \mbox{}^1d\mbox{}^1x)$ \\
for each $x\in {\mathcal C}_j$ and $\mbox{}^1x\in \mbox{}^1{\mathcal
C}_i$.
\par For example, we consider $G$-graded $B$-complexes $({\mathcal C},d)$ and
$(\mbox{}^1{\mathcal C},\mbox{}^1d)$ like $({\mathcal K}(B,X), d)$
and $({\mathcal K}(B,\mbox{}^1X), \mbox{}^1d)$ (see Remark 2.28). By
virtue of Proposition 2.29 Identities $(3.1.1)$, $(3.1.2)$ and
$(3.2.1)$ are satisfied for $D$ naturally induced by $d$ and
$\mbox{}^1d$, since $t_3(g_1,g_2,g_3)\in {\bf \Psi }\subset {\sf
C}(G)\subset N(G)$, $t_n(g_1,...,g_n;u(n),w(n))\in {\bf \Psi }$,
$(g_1g_2)y=g_1(g_2y)$, $(g_1y)g_2=g_1(yg_2)$ and
$y(g_1g_2)=(yg_1)g_2$ for each $g_1$,...,$g_n$ in $G$, $y\in N(G)$,
vectors $u(n)$ and $w(n)$ indicating an order of multiplications
(see Proposition 2.27). \par Thus this example justifies Definition
3.1.
\par In particular, if $\mbox{}^1{\mathcal C}$ has $\mbox{}^1{\mathcal C}_0=X$ and
$\mbox{}^1{\mathcal C}_n=(0)$ for each $n\ne 0$, then $({\mathcal
C}\bigotimes _B\mbox{}^1{\mathcal C})_n={\mathcal C}_n\bigotimes
_BX$ for each $n\in {\bf Z}$, $D=d\otimes I_X$. Therefore ${\mathcal
C}\bigotimes _B B_s$ is isomorphic with ${\mathcal C}$, where $B_s$
is the algebra $B$ considered as the $G$-graded left $B$-module. On
the other side, if ${\mathcal C}_0=P$ and ${\mathcal C}_n=(0)$ for
each $n\ne 0$, then $({\mathcal C}\bigotimes _B \mbox{}^1{\mathcal
C})_n=P\bigotimes _B \mbox{}^1{\mathcal C}_n$ for each $n\in {\bf
Z}$ and $D=I_P\otimes d$.

\par {\bf Proposition 3.3.} {\it Let $({\mathcal C},d)$ and $(\mbox{}^1{\mathcal C},\mbox{}^1d)$ be
$G$-graded $B$-complexes and let ${\mathcal C}$ be a $B$-bimodule
and $\mbox{}^1{\mathcal C}$ be a left $B$-module (or ${\mathcal C}$
be a right $B$-module and $\mbox{}^1{\mathcal C}$ be a $B$-bimodule)
(see Definition 3.1). Then there exists a ${\mathcal T}$-linear
${\bf Z}$-graded map $\hat{h}=\hat{h}({\mathcal
C},\mbox{}^1{\mathcal C})$ of degree $0$ from $H({\mathcal
C})\bigotimes _BH(\mbox{}^1{\mathcal C})$ into $H({\mathcal
C}\bigotimes _B\mbox{}^1{\mathcal C})$.}
\par {\bf Proof.} Consider $x\in {\mathcal Z}_j({\mathcal C})$, $y\in {\mathcal Z}_l(\mbox{}^1{\mathcal C})$, where $j$ and
$l$ are integers (see Definition 3.1 and Remark 3.2). Then the
element $x\otimes y$ belongs to ${\mathcal Z}_{j+l}({\mathcal
C}\bigotimes _B \mbox{}^1{\mathcal C})$ according to Formula
$(3.1.1)$. Therefore $(x+dw)\otimes (y+dz)=x\otimes y+D(w\otimes
y+(-1)^j(x+dw)\otimes z)$ for each $w\in {\mathcal C}_{j+1}$ and
$z\in \mbox{}^1{\mathcal C}_{l+1}$. This induces a so called
canonical ${\mathcal T}$-linear map $\hat{h}_{j,l}: H_j({\mathcal
C})\bigotimes _BH_l(\mbox{}^1{\mathcal C})\to H_{j+l}({\mathcal
C}\bigotimes _B\mbox{}^1{\mathcal C})$ such that
\par $(3.3.1.)$ $\hat{h}_{j,l}(b_ga_u,c_v)={\sf
t}_3(g,u,v)\hat{h}_{j,l}(b_g,a_uc_v)$ and
\par $\hat{h}_{j,l}(a_ub_g,c_v)={\sf
t}_3(u,g,v)a_u\hat{h}_{j,l}(b_g,c_v)$  \par (or
$\hat{h}_{j,l}(b_g,c_va_u)={\sf
t}_3(g,v,u)\hat{h}_{j,l}(b_g,c_v)a_u$ respectively)
\\ for each $a_g\in (H_j({\mathcal C}))_g$, $c_v\in (H_l(\mbox{}^1{\mathcal C}))_v$,
$a_u\in B_u$, $g$, $u$ and $v$ in $G$. The left $B$-bimodule
$H({\mathcal C})\bigotimes _BH(\mbox{}^1{\mathcal C})$ (or right
respectively) is supplied with ${\bf Z}$-gradation such that
\par $(3.3.2)$ $(H({\mathcal C})\bigotimes _BH(\mbox{}^1{\mathcal C}))_n=\sum_{j+l=n}
H_j({\mathcal C})\bigotimes _BH_{l}(\mbox{}^1{\mathcal C})$.
\par Naturally each $H_j({\mathcal C})$ is the (smashly) $G^{|j|+2}$-graded $B$-bimodule
and $H_l(\mbox{}^1{\mathcal C})$ is the (smashly) $G^{|l|+2}$-graded
left $B$-module (or the right $B$-module and the $B$-bimodule
respectively), consequently, $(H({\mathcal C})\bigotimes
_BH(\mbox{}^1{\mathcal C}))_n$ is the $G^{|n|+2}$-graded left
$B$-module (or right respectively) by Lemma 2.5. Therefore, the
family of maps $\hat{h}_{j,l}$ induces a ${\mathcal T}$-linear ${\bf
Z}$-graded map $\hat{h}=\hat{h}({\mathcal C},\mbox{}^1{\mathcal C})$
of degree $0$ from $H({\mathcal C})\bigotimes _BH(\mbox{}^1{\mathcal
C})$ into $H({\mathcal C}\bigotimes _B\mbox{}^1{\mathcal C})$.
\par {\bf Corollary 3.4.} {\it If the conditions of Proposition 3.3 are satisfied and
the $G$-graded $B$-complexes ${\mathcal C}$ and $\mbox{}^1{\mathcal
C}$ are zero from the right, then the $G$-graded $B$-complex
${\mathcal C}\bigotimes _B\mbox{}^1{\mathcal C}$ is zero from the
right and $\hat{h}_{0,0}({\mathcal C},\mbox{}^1{\mathcal C}):
H_0({\mathcal C})\bigotimes _BH_0(\mbox{}^1{\mathcal C})\to
H_0({\mathcal C}\bigotimes _B\mbox{}^1{\mathcal C})$ is bijective.}
\par {\bf Remark 3.5.} Assume that $f: ({\mathcal C},d)\to (\mbox{}^1{\mathcal C},
\mbox{}^1d)$ and $p: (\mbox{}^2{\mathcal C},\mbox{}^2d)\to
(\mbox{}^3{\mathcal C}, \mbox{}^3d)$ are homomorphisms of $G$-graded
$B$-complexes, where the pairs $(({\mathcal C},d),
(\mbox{}^2{\mathcal C}, \mbox{}^2d))$ and $((\mbox{}^1{\mathcal
C},\mbox{}^1d),(\mbox{}^3{\mathcal C}, \mbox{}^3d))$ satisfy the
conditions of Proposition 3.3. Then they induce a homomorphism of
$G$-graded $B$-modules $f\otimes p: {\mathcal C}\bigotimes
_B\mbox{}^2{\mathcal C}\to \mbox{}^1{\mathcal C}\bigotimes _B
\mbox{}^3{\mathcal C}$ such that $(f\otimes p)_n: ({\mathcal
C}\bigotimes _B\mbox{}^2{\mathcal C})_n\to (\mbox{}^1{\mathcal
C}\bigotimes _B \mbox{}^3{\mathcal C})_n$ for each $n\in {\bf Z}$
for $G^{|n|+2}$-graded $B$-modules. This means that the homomorphism
$f\otimes p$ is ${\bf Z}$-graded of zero degree. For derivations $D$
and $\mbox{}^1D$ of ${\mathcal C}\bigotimes _B\mbox{}^2{\mathcal C}$
and $\mbox{}^1{\mathcal C}\bigotimes _B \mbox{}^3{\mathcal C}$
respectively we get that
\par $(f\otimes p)(D(x\otimes y))=f(dx)\otimes
p(y)+(-1)^jf(x)\otimes p(\mbox{}^2dy)=\mbox{}^1df(x)\otimes
p(y)+(-1)^jf(x)\otimes \mbox{}^3dp(y)=\mbox{}^1D(f(x)\otimes p(y))$
for each $x\in {\mathcal C}_j$, $y\in \mbox{}^2{\mathcal C}_l$, $j$
and $l$ in ${\bf Z}$. In view of Proposition 3.3 this provides a
commutative diagram
\par  $(3.5.1) \quad
$ $~ ~ {H({\mathcal C})\bigotimes _ BH(\mbox{}^2{\mathcal
C})}_{\overrightarrow{\hspace{0.8cm}\hat{h}({\mathcal
C},\mbox{}^2{\mathcal C}) \hspace{0.8cm}}} H({\mathcal C}\bigotimes
_B\mbox{}^2{\mathcal C})$
\par \hspace{1.0cm} $H(f)\otimes H(p)\downarrow \hspace{5.5cm} \downarrow H(f\otimes p)$
\par \hspace{1.5cm} $~  {H(\mbox{}^1{\mathcal C})\bigotimes _BH(\mbox{}^3{\mathcal C})}_{
\overrightarrow{\hspace{0.8cm}\hat{h}(\mbox{}^1{\mathcal
C},\mbox{}^3{\mathcal C})\hspace{0.8cm}}} {H(\mbox{}^1{\mathcal
C}\bigotimes _B\mbox{}^3{\mathcal C})}$.

\par {\bf Proposition 3.6.} {\it Let ${\mathcal C}$, $\mbox{}^1{\mathcal C}$, $\mbox{}^2{\mathcal C}$ and
$\mbox{}^3{\mathcal C}$ be $G$-graded $B$-complexes, let the pairs
$({\mathcal C},\mbox{}^2{\mathcal C})$ and $(\mbox{}^1{\mathcal C},
\mbox{}^3{\mathcal C})$ satisfy the conditions of Proposition 3.3
and let $f: {\mathcal C}\to \mbox{}^1{\mathcal C}$, $\mbox{}^1f:
{\mathcal C}\to \mbox{}^1{\mathcal C}$, $p: \mbox{}^2{\mathcal C}\to
\mbox{}^3{\mathcal C}$, $\mbox{}^1p: \mbox{}^2{\mathcal C}\to
\mbox{}^3{\mathcal C}$ be $B$-generic homomorphisms of these
complexes. Then two homomorphisms $f\otimes p$ and
$\mbox{}^1f\otimes \mbox{}^1p$ from ${\mathcal C}\bigotimes
_B\mbox{}^2{\mathcal C}$ to $\mbox{}^1{\mathcal C}\bigotimes
_B\mbox{}^3{\mathcal C}$ are $B$-generic.
\par $(3.6.1).$ If $f$ and $p$ are homotopic to $\mbox{}^1f$ and
$\mbox{}^1p$ respectively, then two homomorphisms $f\otimes p$ and
$\mbox{}^1f\otimes \mbox{}^1p$ are homotopic.
\par $(3.6.2).$ If $f$ and $p$ are homotopisms, then $f\otimes p$ is
a homotopism.
\par  $(3.6.3).$ If either ${\mathcal C}$ or $\mbox{}^2{\mathcal C}$ is homotopic to zero, then
${\mathcal C}\bigotimes _B\mbox{}^2{\mathcal C}$ is homotopic to
zero.}
\par {\bf Proof.} From $f_{\bf \iota }(B)=B$, $p_{\bf \iota }(B)=B$, $\mbox{}^1f(B)=B$,
$\mbox{}^1p(B)=B$ provided by the conditions of this proposition it
follows that $(f_{\bf \iota }\otimes p_{\bf \iota })(B\otimes
_BB)=B\otimes _BB$ and $(\mbox{}^1f_{\bf \iota }\otimes
\mbox{}^1p_{\bf \iota })(B\otimes _BB)=B\otimes _BB$, consequently,
the homomorphisms $f\otimes p$ and $\mbox{}^1f\otimes \mbox{}^1p$
are $B$-generic.
\par If $f$ and $\mbox{}^1f$ are homotopic
to $p$ and $\mbox{}^1p$ respectively, then there exist ${\bf
Z}$-graded $B$-generic homomorphisms $s: {\mathcal C}\to
\mbox{}^1{\mathcal C}$ and $\mbox{}^1s: \mbox{}^2{\mathcal C}\to
\mbox{}^3{\mathcal C}$ of degree $1$ such that $f-\mbox{}^1f=ds+sd$
and $p-\mbox{}^1p=d\mbox{ }^1s+\mbox{}^1sd$, where derivations of
the $G$-graded $B$-complexes ${\mathcal C}$, $\mbox{}^1{\mathcal
C}$, $\mbox{}^2{\mathcal C}$ and $\mbox{}^3{\mathcal C}$ are shortly
denoted by $d$. Therefore, there exists a ${\bf Z}$-graded
homomorphism $S: {\mathcal C}\bigotimes _B\mbox{}^2{\mathcal C}\to
\mbox{}^1{\mathcal C}\bigotimes _B\mbox{}^3{\mathcal C}$ of degree
$1$ such that $S(x\otimes y)=s(x)\otimes p(y)+(-1)^j
~\mbox{}^1f(x)\otimes \mbox{}^1s(y)$ for each $x\in {\mathcal C}_j$,
$y\in \mbox{}^2{\mathcal C}_l$, $j$ and $l$ in ${\bf Z}$. Since
$\mbox{}^1f$, $p$, $s$, $\mbox{}^1s$ are $B$-generic, then $S$ is
$B$-generic. For derivations $D$ of the $G$-graded $B$-complexes
${\mathcal C}\bigotimes _B\mbox{}^2{\mathcal C}$ and
$\mbox{}^1{\mathcal C}\bigotimes _B\mbox{}^3{\mathcal C}$ this
gives: $(DS+SD)(x\otimes y)= f(x)\otimes p(y)-\mbox{}^1f(x)\otimes
\mbox{}^1p(y)$ for each $x\in {\mathcal C}_j$, $y\in
\mbox{}^2{\mathcal C}_l$, $j$ and $l$ in ${\bf Z}$. Thus
$DS+SD=f\otimes p-\mbox{}^1f\otimes \mbox{}^1p$, that means that two
homomorphisms $f\otimes p$ and $\mbox{}^1f\otimes \mbox{}^1p$ are
homotopic.
\par If $f$ and $p$ are homotopisms, then there exist $B$-generic
homomorphisms of complexes $\eta : \mbox{}^1{\mathcal C}\to
{\mathcal C}$ and $\mbox{}^1\eta : \mbox{}^3{\mathcal C}\to
\mbox{}^2{\mathcal C}$ such that $f\circ \eta $, $\eta \circ f$,
$p\circ \mbox{}^1\eta $, $\mbox{}^1\eta \circ p$ are homotopic to
$id_{\mbox{}^1{\mathcal C}}$, $id_{{\mathcal C}}$,
$id_{\mbox{}^3{\mathcal C}}$, $id_{\mbox{}^2{\mathcal C}}$
respectively. From $(3.6.1)$ it follows that $(f\otimes p)\circ
(\eta \otimes \mbox{}^1\eta )$ is homotopic to
$id_{\mbox{}^1{\mathcal C}}\otimes id_{\mbox{}^3{\mathcal
C}}=id_{\mbox{}^1{\mathcal C}\otimes \mbox{}^3{\mathcal C}}$, while
$(\eta \otimes \mbox{}^1\eta )\circ (f\otimes p)$ is homotopic to
$id_{{\mathcal C}}\otimes id_{\mbox{}^2{\mathcal C}}=id_{{\mathcal
C}\otimes \mbox{}^2{\mathcal C}}$. This implies that $f\otimes p$ is
the homotopism.
\par The last assertion $(3.6.3)$ of this proposition follows from
$(3.6.2)$ in particular for either $\mbox{}^2{\mathcal C}=0$ or
$\mbox{}^3{\mathcal C}=0$ respectively.

\par \textbf{Conclusion 3.7.} The results of this article
can be used for further studies of cohomology theory of
nonassociative algebras and noncommutative manifolds with metagroup
relations. Then it is interesting to mention possible applications
in mathematical coding theory, analysis of information flows and
their technical realizations
\cite{blautrctb,maglebrtj19,srwseabm14}, because frequently codes
are based on binary systems and algebras. Indeed, metagroup
relations are weaker, than relations in groups. Therefore, a code
complexity can increase by using nonassociative algebras with
metagroup relations in comparison with group algebras or Lie
algebras. \par Besides applications of cohomologies outlined in the
introduction they also can be used in mathematical physics and
quantum field theory. They can be used for further studies of
nonassociative algebra cohomologies, structure of nonassociative
algebras, operator theory and spectral theory over Cayley-Dickson
algebras, PDEs, noncommutative analysis, noncommutative geometry,
mathematical physics, their applications in the sciences. This also
can be applied to cohomologies of PDEs and solutions of PDEs with
boundary conditions which can have a practical importance
\cite{pommb}.

\section{Appendix. Nonassociative algebras with metagroup relations
and their modules.}

\par In this appendix definitions from the previous article
\cite{ludkctnaax19} are reminded. A reader familiar with it can skip
this appendix.

\par {\bf Definition A1.}
Let $G$ be a set with a single-valued binary operation
(multiplication)  $G^2\ni (a,b)\mapsto ab \in G$ defined on $G$
satisfying the conditions: \par $(A1.1)$ for each $a$ and $b$ in $G$
there is a unique $x\in G$ with $ax=b$ and \par $(A1.2)$ a unique
$y\in G$ exists satisfying $ya=b$, which are denoted by \par
$x=a\setminus b=Div_l(a,b)$ and $y=b/a=Div_r(a,b)$ correspondingly,
\par $(A1.3)$ there exists a neutral (i.e. unit) element $e_G=e\in G$:
\par $~eg=ge=g$ for each $g\in G$.
\par The set of all elements $h\in G$
commuting and associating with $G$:
\par $(A1.4)$ $Com (G) := \{ a\in G: \forall b\in G, ~ ab=ba \} $,
\par $(A1.5)$ $N_l(G) := \{a\in G: \forall b\in G, \forall c\in G, ~ (ab)c=a(bc) \}
$,
\par $(A1.6)$ $N_m(G) := \{a\in G: \forall b\in G, \forall c\in G, ~ (ba)c=b(ac)
\} $,
\par $(A1.7)$ $N_r(G) := \{a\in G: \forall b\in G, \forall c\in G, ~ (bc)a=b(ca)
\} $,
\par $(A1.8)$ $N(G) := N_l(G)\cap N_m(G)\cap N_r(G)$; \par ${\mathcal C}(G) := Com (G)\cap N(G)$
is called the center ${\mathcal C}(G)$ of $G$.
\par We call $G$ a metagroup if a set $G$ possesses a single-valued binary operation
and satisfies conditions $(A1.1)$-$(A1.3)$ and
\par $(A1.9)$ $(ab)c={\sf t}_3(a,b,c)a(bc)$ \\ for each
$a$, $b$ and $c$ in $G$, where ${\sf t}_3(a,b,c)\in {\bf \Psi }$, $
~ {\bf \Psi }\subset {\mathcal C}(G)$; \\ where ${\sf t}_3$ shortens
a notation ${\sf t}_{3,G}$, where ${\bf \Psi }$ denotes a (proper or
improper) subgroup of ${\mathcal C}(G)$.
\par Then $G$ will be called a central metagroup if in addition to $(A1.9)$ it
satisfies the condition:
\par $(A1.10)$ $ab={\sf t}_2(a,b)ba$ \\ for each $a$ and $b$ in $G$,
where ${\sf t}_2(a,b)\in {\bf \Psi }$.
\par Particularly, $Inv_l(a)=Div_l(a,e)$ is a left inversion,
$Inv_r(a)=Div_r(a,e)$ is a right inversion.
\par In view of the nonassociativity of $G$ in general
a product of several elements of $G$ is specified as usually by
opening "$($" and closing "$)$" parentheses. For elements
$a_1$,...,$a_n$ in $G$ we shall denote shortly by $ \{ a_1,...,a_n
\}_{q(n)}$ the product, where a vector $q(n)$ indicates an order of
pairwise multiplications of elements in the row $a_1,...,a_n$ in
braces in the following manner. Enumerate positions: before $a_1$ by
$1$, between $a_1$ and $a_2$ by $2$,..., by $n$ between $a_{n-1}$
and $a_n$, by $n+1$ after $a_n$. Then put $q_j(n)=(k,m)$ if there
are $k$ opening "$($" and $m$ closing "$)$" parentheses in the
ordered product at the $j$-th position of the type $)...)(...($,
where $k$ and $m$ are nonnegative integers, $q(n)
=(q_1(n),....,q_{n+1}(n))$ with $q_1(n)=(k,0)$ and
$q_{n+1}(n)=(0,m)$.
\par As traditionally $S_n$ denotes the symmetric group of the set $ \{ 1,
2,..., n \} $. Henceforth, maps and functions on metagroups are
supposed to be single-valued if something other will not be
specified.
\par Let $\psi : G\to G$ be a bijective surjective map satisfying
the following condition: $\psi (ab)=\psi (a)\psi (b)$ for each $a$
and $b$ in $G$. Then $\psi $ is called an automorphism of the
metagroup $G$.

\par {\bf Definition A2.} Let $A$ be an algebra over an
associative unital ring $\mathcal T$ such that $A$ has a natural
structure of a $({\mathcal T}, {\mathcal T})$-bimodule with a
multiplication map $A\times A\to A$, which is right and left
distributive $a(b+c)=ab+ac$, $(b+c)a=ba+ca$, also satisfying the
following identities $r(ab)=(ra)b$, $~(ar)b=a(rb)$, $~(ab)r=a(br)$,
$~s(ra)=(sr)a$ and $(ar)s=a(rs)$ for any $a$, $b$ and $c$ in $A$,
$r$ and $s$ in ${\mathcal T}$. Let $G$ be a metagroup and $\mathcal
T$ be an associative unital ring.
\par Henceforth the ring ${\mathcal T}$ will be supposed commutative, if
something other will not be specified.
\par Then by ${\mathcal
T}[G]$ is denoted a metagroup algebra over $\mathcal T$ of all
formal sums $s_1a_1+...+s_na_n$ satisfying conditions $(1-3)$ below,
where $n$ is a positive integer, $s_1$,...,$s_n$ are in $\mathcal T$
and $a_1$,...,$a_n$ belong to $G$:
\par $(A2.1)$ $sa=as$ for each $s$ in $\mathcal T$ and $a$ in $G$,
\par $(A2.2)$ $s(ra)=(sr)a$ for each $s$ and $r$ in $\mathcal T$, and
$a\in G$, \par $(A2.3)$ $r(ab)=(ra)b$, $~(ar)b=a(rb)$,
$~(ab)r=a(br)$ for each $a$ and $b$ in $G$, $r\in {\mathcal T}$.

\par {\bf Definition A3.} Let $\mathcal R$ be a ring, which may be nonassociative relative to the
multiplication. If there exists a mapping ${\mathcal R}\times M\to
M$, $~{\mathcal R}\times M\ni (a,m)\mapsto am\in M$ such that
$a(m+k)=am+ak$ and $(a+b)m=am+bm$ for each $a$ and $b$ in ${\mathcal
R}$, $m$ and $k$ in $M$, then $M$ will be called a generalized left
${\mathcal R}$-module or shortly: left ${\mathcal R}$-module or left
module over $\mathcal R$.
\par If ${\mathcal R}$ is a unital ring and $1m=m$ for each $m\in M$,
then $M$ is called a left unital module over $\mathcal R$, where $1$
denotes the unit element in the ring ${\mathcal R}$. Symmetrically
is defined a right ${\mathcal R}$-module. \par If $M$ is a left and
right ${\mathcal R}$-module, then it is called a two-sided
${\mathcal R}$-module or a $({\mathcal R},{\mathcal R})$-bimodule.
If $M$ is a left ${\mathcal R}$-module and a right ${\mathcal
S}$-module, then it is called a $({\mathcal R}, {\mathcal
S})$-bimodule.
\par A two-sided module $M$ over $\mathcal R$ is called cyclic, if an
element $y\in M$ exists such that $M={\mathcal R}(y{\mathcal R})$
and $M=({\mathcal R}y){\mathcal R}$, where ${\mathcal R}(y{\mathcal
R}) = \{ s(yp): ~ s, p\in {\mathcal R} \} $ and $({\mathcal
R}y){\mathcal R} = \{ (sy)p: ~ s, p\in {\mathcal R} \} $.
\par Let $G$ be a metagroup. Take a metagroup algebra $A={\mathcal T}[G]$ and a two-sided $A$-module $M$, where
$\mathcal T$ is an associative unital ring (see Definition 2). Let
$M_g$ be a two-sided ${\mathcal T}$-module for each $g\in G$, where
$G$ is the metagroup. Let $M$ have the decomposition $M=\sum_{g\in
G}M_g$ as a two-sided ${\mathcal T}$-module. Let also $M$ satisfy
the following conditions:
\par $(A3.1)$ $hM_g=M_{hg}$ and $M_gh=M_{gh}$,
\par $(A3.2)$ $(bh)x_g=b(hx_g)$ and $x_g(bh)=(x_gh)b$ and $bx_g=x_gb$,
\par $(A3.3)$ $(hs)x_g={\sf t}_3(h,s,g) h(sx_g)$ and $(hx_g)s= {\sf t}_3(h,g,s) h(x_gs)$
and \par $(x_gh)s={\sf t}_3(g,h,s) x_g(hs)$ \\ for every $h, g, s$
in $G$ and $b\in \mathcal T$ and $x_g\in M_g$. Then a two-sided
$A$-module $M$ satisfying conditions $(A3.1)$-$(A3.3)$ will be
called smashly $G$-graded. For short it also will be said
"$G$-graded" instead of "smashly $G$-graded". In particular, if the
module $M$ is $G$-graded and splits into a direct sum
$M=\bigoplus_{g\in G}M_g$ of two-sided ${\mathcal T}$-submodules
$M_g$, then we will say that that $M$ is directly $G$-graded. For a
nontrivial (nonzero) $G$-graded module $X$ with the nontrivial
metagroup $G$ it will be supposed that there exists $g\in G$ such
that $X_g\ne X_e$, if something other will not be outlined. \par
Similarly are defined $G$-graded left and right $A$-modules.
Henceforward, speaking about $A$-modules (left, right or two-sided)
it will be supposed that they are $G$-graded and it will be written
for short "an $A$-module" instead of "a $G$-graded $A$-module",
unless otherwise specified.
\par If $P$ and $N$ are left
$A$-modules and a homomorphism $\gamma : P\to N$ is such that
$\gamma (ax)=a\gamma (x)$ for each $a\in A$ and $x\in P$, then
$\gamma $ is called a left $A$-homomorphism. Analogously are defined
right $A$-homomorphisms for right $A$-modules. For two-sided $A$
modules a left and right $A$-homomorphism is called an
$A$-homomorphism. \par For left $\mathcal T$-modules $M$ and $N$ by
$Hom_{\mathcal T}(M,N)$ is denoted a family of all left $\mathcal
T$-homomorphisms from $M$ into $N$. A similar notation is used for a
family of all $\mathcal T$-homomorphisms (or right $\mathcal
T$-homomorphisms) of two-sided $\mathcal T$-modules (or right
$\mathcal T$-modules correspondingly). If an algebra $A$ is
specified it may be written shortly a homomorphism instead of an
$A$-homomorphism.

\par \textbf{Conflicts of Interest:} The author declares no
conflict of interest.

\end{document}